\documentclass[12pt]{article}
\usepackage{ amsfonts,
amsmath, amssymb,amsgen, amsthm, amscd}

\def\GL{\mathop{\rm GL}\nolimits}

\def\Im{\mathop{\rm Im}\nolimits}
\def\Re{\mathop{\rm Re}\nolimits}
\def\det{\mathop{\rm det}\nolimits}

 \def\Hom{\mathop{\rm Hom}\nolimits}

\def\Log{\mathop{\rm Log}\nolimits}
\def\log{\mathop{\rm log}\nolimits}

\def\sign{\mathop{\rm sign}\nolimits}

\def\SL{\mathop{\rm SL}\nolimits}
\def\PSL{\mathop{\rm PSL}\nolimits}

\def\Ac{{\cal A}}

\def\Ec{{\cal E}}

\def\Gc{{\cal G}}
\def\Hc{{\cal H}}

\def\Lc{{\cal L}}
\def\Mc{{\cal M}}

\def\Rc{{\cal R}}
\def\Sc{{\cal S}}

\def\a{\alpha}
\def\b{\beta}
\def\d{\delta}
\def\D{\Delta}
\def\g{\gamma}
\def\e{\epsilon}
\def\G{\Gamma}
\def\lb{\lambda}

\def\om{\omega}
\def\Om{\Omega}

\def\s{\sigma}

\def\t{\tau}
\def\th{\theta}
\def\ve{\varepsilon}

\def\ab{\bf a}
\def\eb{\bf e}
\def\kb{\bf k}
\def\mb{\bf m}
\def\xb{\bf x}
\def\yb{\bf y}

\def\Bb{\bf B}

\def\0b{\bf 0}

\def\fl{\forall}
\def\ify{\infty}

\def\ot{\otimes}

\def\ra{\rightarrow}

\def\sbs{\subset}

\def\wt{\widetilde}

\def\Qb{{\mathbb Q}}
\def\Cb{{\mathbb C}}

\def\Hb{{\mathbb H}}
\def\Nb{{\mathbb N}}
\def\Rb{{\mathbb R}}
\def\Tb{{\mathbb T}}
\def\Zb{{\mathbb Z}}

\newtheorem{theorem}{Theorem}
\newtheorem{remark}[theorem]{Remark}
\newtheorem{proposition}[theorem]{Proposition}
\newtheorem{lemma}[theorem]{Lemma}

\def\build#1_#2^#3{\mathrel{
\mathop{\kern 0pt#1}\limits_{#2}^{#3}}}

\numberwithin{equation}{section}
\begin{document}

\title{{\bf Transgressions of the Godbillon-Vey class
and Rademacher functions}}
\author{Alain Connes \\
        Coll\`ege de France \\
        3 rue d'Ulm \\
        75005 Paris, France \\
\and
        Henri Moscovici\thanks{Research
    supported by the National Science Foundation
    award no. DMS-0245481.} \\
    Department of Mathematics \\
    The Ohio State University \\
    Columbus, OH 43210, USA
    }

\date{ \ }

\maketitle

\section*{Introduction}

In earlier work \cite{CMmha, CMrc}
we investigated a surprising interconnection between the
transverse geometry of codimension $1$ foliations and modular forms.
At the core of this interplay lies the Hopf algebra $\Hc_{1}$,
the first in a series
of Hopf algebras $\Hc_n$ that were found \cite{CMhti}
to determine the affine transverse geometry
of codimension $n$ foliations.
The periodic Hopf-cyclic cohomology of
$\Hc_{1}$ is generated by two classes, $[\d_1] $ for the
odd component and $[RC_1]$ for the even
component. The tautological action of $\Hc_{1}$ on the \'etale groupoid algebra
$\Ac_{\Gc}$ associated to the frame bundle of a codimension $1$ foliation
preserves (up to a character) the canonical trace on $\Ac_{\Gc}$, and thus
gives rise to a characteristic
homomorphism in cyclic cohomology. This homomorphism maps $[\d_1] $
to the Godbillon-Vey class and $[RC_1]$ to the transverse fundamental class.
\medskip

The starting point of the investigation in \cite{CMmha} was the realization that
the Hopf algebra $\Hc_{1}$ can be made
to act on the crossed product $\Ac_{\Qb}$ of the algebra of modular forms of
all levels by $\GL^+ (2,{\mathbb Q})$, via a natural connection provided by
the Ramanujan operator on modular forms,
thus conferring a symmetry structure
to the space of lattices modulo the action of the Hecke correspondences.
Although the algebra $\Ac_{\Qb}$ no longer has an invariant trace,
we used an ad hoc pairing with modular symbols to
convert the Godbillon-Vey class $[\d_1] \in
HC^1 (\Hc_{1})$ into the Euler class of $\GL^+ (2,{\mathbb Q})$.
\medskip

In this paper
we provide a completely conceptual explanation for the above pairing
and at the same time extend it to the higher weight case. This is
achieved by constructing out of modular symbols
$\Hc_{1}$-invariant $1$-traces
that support characteristic maps for certain actions
of $\Hc_{1}$ on $\Ac_{\Qb}$, canonically associated to modular forms.
Moreover, we
show that the image of the Godbillon-Vey
class through these characteristic homomorphisms, obtained by the
cup product  between  $[\d_1] $ and the invariant $1$-traces,
transgresses to secondary data.
For the projective action determined by the Ramanujan connection
the transgression takes place within the Euler class, in a manner that ressembles
the $K$-homological transgression in the context of $SU_q(2)$ \cite{C2},
and leads to the classical
Rademacher function \cite{RG}. For the actions associated to cusp forms
of higher weight
the transgressed classes implement the Eichler-Shimura isomorphism.
The actions corresponding to
Eisenstein series give rise by transgression to higher Dedekind sums
and generalized Rademacher functions (\cite{Ste, Nak}), or equivalently
to the Eisenstein cocycle of \cite{Ste}.
\medskip

Generalized Dedekind sums have been related to
special values of $L$-functions in
the work of C. Meyer on the class-number formula \cite{Me1, Me2},
and higher Dedekind sums appear
in the work of Siegel  \cite{Si1, Si2}
and Zagier  \cite{Za1} on the values
at non-positive integers of partial zeta functions over
real quadratic fields.
Eisenstein cocycles were employed by Stevens \cite{Ste}
and by Sczech \cite{Sc} in order
to compute these values more efficiently.
The fact that these notions can be interpreted as secondary invariants
is reminiscent of the secondary nature of the
regulator invariants (cf. e.g.~\cite{DHZ}) that are involved
in the expression of the special values at non-critical
points of $L$-functions associated to number fields
(see e.g.~\cite{KZ}, \cite{Za2}).
 \bigskip

 \tableofcontents

\section{The standard modular Hopf action}

In this preliminary section, we briefly review the basic facts
(cf. \cite{CMhti}, \cite{CMmha}) concerning the
Hopf algebra $\Hc_1$ and its standard Hopf action on
the crossed product $\Ac_{\Qb}$ of the algebra of modular forms of
all levels by $\GL^+ (2,{\mathbb Q})$,
associated to the Ramanujan connection.

\subsection{The Hopf algebra $\Hc_1$ and its cyclic classes}
 We start by recalling the definition of the Hopf
algebra $\Hc_1$. As an algebra, it coincides with
the universal enveloping algebra of the Lie algebra with basis $\{
X,Y,\d_n \, ; n \geq 1 \}$ and brackets
\begin{equation*} \label{pres}
[Y,X] = X \, , \, [Y , \d_n ] = n \, \d_n \, , \, [X,\d_n] =
\d_{n+1} \, , \, [\d_k , \d_{\ell}] = 0 \, , \quad n , k , \ell
\geq 1 \, .
\end{equation*}
As a Hopf algebra, the coproduct $\, \D : \Hc_1 \ra \Hc_1 \ot
\Hc_1 \,$ is determined by
\begin{eqnarray*}
\D \,  Y \, &=& \, Y \ot 1 + 1 \ot Y \, , \qquad \qquad  \D \,  \d_1 \, = \, \d_1 \ot 1 + 1 \ot \d_1
 \nonumber \\
 \D \,  X \, &=& \, X \ot 1 + 1 \ot X + \d_1 \ot Y
\end{eqnarray*}
and the multiplicativity property
\begin{equation*}
\D (h^1 \, h^2) = \D h^1 \cdot \D h^2 \, , \quad h^1 , h^2 \in
\Hc_1 \, \, ;
\end{equation*}
the antipode is determined by
\begin{equation*}
S(Y) = -Y , \quad S(X) = -X + \d_1 Y  , \quad S(\d_1) = - \d_1
\end{equation*}
and the anti-isomorphism property
\begin{equation*}
S (h^1 \, h^2) = S(h^2) \, S (h^1) \, , \quad h^1 , h^2 \in \Hc_1
\, \, ;
\end{equation*}
finally, the counit is
\begin{equation*}
\ve (h) = \hbox{constant term of} \quad h \in \Hc_1 \, .
\end{equation*}
\bigskip

The modular character $\,\d \in \Hc_1^* \,$, determined by
$$\d (Y) = 1, \quad \d (X) = 0 , \quad \d (\d_n) = 0 \, ,
$$
together with the unit of $1 \in \Hc_1$ forms a
 modular pair in involution $(\d, 1)$, and thus
the Hopf-cyclic cohomology $HC_{(\d, 1)}^1 \, (\Hc_1)$  is well-defined
(for definitions, see \cite{CMhti, CMcchs}).

\noindent The element $\, \d_{1} \in \Hc_{1} \,$ is a Hopf-cyclic cocycle,
which gives a nontrivial class
$$
 [\d_{1}] \, \in \, HC_{(\d, 1)}^1 \, (\Hc_1) \,
$$
in the Hopf-cyclic cohomology of $\Hc_1$
with respect to the modular pair $( \d, 1)$.
Its periodic image generates the periodic group
$ \,HP^{1} \, (\Hc_1; \d, 1) \,$,
and represents the universal Godbillon-Vey class
(cf. \cite[Prop. 3]{CMmha}).
\smallskip

The even component of the
periodic cyclic cohomology group $HP^{0}(\Hc_1 ; \d, 1)$
is generated by the ``transverse
fundamental class'', represented by the Hopf-cyclic $2$-cocycle
\begin{equation*} \label{fund}
RC_1 \,  := \, X \ot Y - Y \ot X - \d_1 \,  Y \ot Y \, .
\end{equation*}
(See \cite{CMrc} for the explanation of the notation.)
\smallskip

There is one other Hopf-cyclic $1$-cocycle,
intimately related to the classical Schwarzian,
which plays a prominent role in the transverse geometry of
modular Hecke algebras (cf. \cite{CMmha, CMrc}). It is given
by the primitive element
$$
\d'_{2} := \d_2- \frac{1}{2}\d_1^2 \, \in \Hc_{1} \, .
$$
Its periodic class vanishes because
$ \d'_2 \, = \, B\, (c)$, where $c$ is the following Hochschild
$2$-cocycle:
$$ c \, := \, \delta_1 \ot X + \frac{1}{2} \, \delta_1^2 \ot Y \, .
$$
 \bigskip

 \subsection{Standard  modular action of $\Hc_1$}
 The notation being as in \cite[\S 1]{CMmha}, we
 form the  crossed product algebra
$$ \Ac_{G^+ (\mathbb{Q})} \, = \, \Mc \ltimes \GL^+ (2,{\Qb})\, ,
$$
where $\Mc$ is the algebra of (holomorphic) modular forms
of all levels. The product of two elements in $ \Ac_{G^+ (\mathbb{Q})} $,
$$
 a^0 = \sum_{\a} f^0_{\a} U_{\a}   \quad \text{and}
 \quad a^1=\sum_{\b} f^1_{\b} U_{\b}   \, ,
 $$
is given by the convolution rule
$$
a^0 \, a^1 \, = \,  \sum_{\a, \b} f^0_{\a} \, f^1_{\b} | \a^{-1} \,U_{\a \b}  \, .
$$
We recall (see  \cite[Prop. 7]{CMmha}) that there is a unique
Hopf action of
the Hopf algebra $\, \Hc_1 \,$ on  $ \Ac_{G^+ (\mathbb{Q})} $
determined by letting the generators $\{ Y, X, \d_1\}$
of $\Hc_1$ act
on monomials $f  U_{\g}^* \in \Ac_{G^+ (\mathbb{Q})} \,$
as follows:
\begin{equation} \label{Y}
Y (f  U_{\g}^*) = Y(f)\, U_{\g}^* , \quad \text{where} \quad
Y(f)  =  \frac{w(f)}{2} \, f ,  \quad w(f) = \text{weight}(f) ;
\end{equation}
\begin{equation}  \label{X}
X (f  U_{\g}^*) =  X(f)\, U_{\g}^* , \quad \text{where} \quad
X  = \frac{1}{2  \pi i} \,
\frac{d}{dz} - \frac{1}{2  \pi i} \, \frac{d}{dz} (\log \eta^4)
\, Y
\end{equation}
and $\eta$ stands for the Dedekind $ \eta$-function,
\begin{equation}  \nonumber
 \eta^{24} (z) \, = \, q \, \prod_{n=1}^{\infty} (1 - q^n)^{24} \, , \qquad q
= e^{2\pi i z} ;
\end{equation}
lastly,
\begin{equation} \label{d1}
\d_1 (f  U_{\g}^*) = \mu_{\g} \, f
U_{\g}^* \, ,
\end{equation}
with the factor $\mu_{\g}$ given by the expression
\begin{equation} \label{mu}
\mu_{\g} \, (z) \, = \, \frac{1}{12  \pi i} \, \frac{d}{dz} \log
\frac{\D | \g}{\D} \, = \, \frac{1}{2  \pi i} \, \frac{d}{dz} \log
\frac{\eta^4 | \g}{\eta^4} \, ;
\end{equation}
equivalently,
\begin{equation} \label{mug}
 \mu_{\g} \, (z) = \frac{1}{2 \, \pi^2}
 \left(  G_2 | \g \, ( z) -  G_2 ( z)
+ \frac{2  \pi i \, c}{cz+d} \right) ,
\end{equation}
where
\begin{equation*} \label{G*}
 G_2 (z) \, = \, 2 \zeta (2) \, + \,
2 \sum_{m \geq 1} \, \sum_{n \in \mathbb{Z}} \frac{1}{(mz + n)^2}
\, = \, \frac{\pi^2}{3} \, - \, 8 \pi^2 \sum_{m, n \geq 1} m e^{2
\pi i m n z}
\end{equation*}
is the quasimodular holomorphic Eisenstein series of weight 2.
The factor $ \mu_{\g} $ can further be expressed as the difference
\begin{equation}  \label{gb0}
 \mu_{\g} \, = \,  2 ( \phi_{\0b} | \g \, - \, \phi_{\0b} ) \, , \quad
 \phi_{\0b}= \frac{1}{4 \pi^2}G_{\0b} \, ,
\end{equation}
where  $G_{\0b}$
is the modular (but nonholomorphic) weight $2$
Eisenstein series
$$
\, G_{\0b} (z) \, = \, G_{\0b} (z, 0)
$$
obtained by taking the value at $\, s = 0 \,$ of the analytic
continuation of the series
\begin{eqnarray}
G_{\0b} (z, s) \, &=& \sum_{(m, n) \in \mathbb{Z}^2 \setminus 0}
(m z + n)^{-2} \, |m z + n|^{-s}  \nonumber \\
  &=& 2 \zeta (2 + s) \, + \,
2 \sum_{m \geq 1} \, \sum_{n \in \mathbb{Z}} (m z + n)^{-2} \, |m
z + n|^{-s} \, , \quad \Re s > 0 \, . \nonumber
\end{eqnarray}
It is related to $\, G_2 \,$ by the identity
\begin{equation*} \label{GG*}
G_2 (z) \, = \, G_{\0b} (z) \, + \, \frac{2  \pi i}{z - \bar{z}} \, .
\end{equation*}
\smallskip

The equation (\ref{gb0}) shows that the range of $\, \mu \,$ is
contained in the  space  $\Ec_2 (\mathbb{Q}) $ of
 weight $2$ Eisenstein series whose constant term in the $q$-expansion
 at each cusp is rational.  We recall (following  \cite[\S 2.4]{St})
 Hecke's construction \cite{He} of a lattice of generators
 for the $\mathbb{Q}$-vector space $\, \Ec_2 (\mathbb{Q})$.

 For $\, {\ab} = (a_1, a_2) \in ( \mathbb Q \slash
\mathbb Z )^2  $ and  $z \in \Hb $ fixed, the series
\begin{equation*}
 G_{\ab} (z, s) := \, \sum_{ {\mb} \neq \0b ,
\, {\mb} \equiv {\ab} \, (\text{mod} \,1)}
 (m_{1} z + m_{2})^{-2} |m_{1} z + m_{2}|^{-s} \, , \quad \Re s > 0 \, \, ;
\end{equation*}
defines a function that can be analytically
continued beyond $\Re s = 0 $, which allows to define
\begin{equation*}
 G_{\ab} (z) := \,  G_{\ab} (z, 0) \, .
\end{equation*}
Furthermore, one has
\begin{equation*}
  G_{\ab} | \g \, = \,  G_{{\ab} \cdot \g} \, , \quad \fl \, \g \in
  \G (1) \, ,
\end{equation*}
which shows that $ G_{\ab} (z) $ behaves like a weight $2$ modular
form of some level $N$. However it is
only quasi-holomorphic, in the sense that the function
$$
z \mapsto G_{\ab} (z) \, + \, \frac{2 \pi i}{z - \bar{z}} 
$$
is holomorphic in $ z \in \Hb$. Moreover,
the difference
\begin{equation*}
     \wp_{\ab} (z) \, = \,  G_{\ab} (z) \, - \,  G_{\0b} (z)
\end{equation*}
is precisely the ${\ab}-$division value of the Weierstrass 
$\wp$-function, and the collection of functions
\begin{equation*}
 \left\{ \wp_{\ab} \, \, ; \,
\, {\ab} \in \left( \frac{1}{N} \mathbb Z \slash \mathbb Z
\right)^2 \setminus {\0b} \right\}
\end{equation*}
generates the space of weight $2$ Eisenstein series of level $N$.
\smallskip

\noindent In order to obtain a set of generators for $\Ec_2 (\mathbb Q)$, one
considers  the additive characters $\, \chi_{\xb} :
 \left( \frac{1}{N} \mathbb Z \slash \mathbb Z \right)^2 \ra
 {\mathbb C}^{\times} $ defined by
\begin{equation} \label{chi}
  \chi_{\xb} \left( \frac{\ab}{N} \right) \, := \,
 e^{2  \pi i\,  (a_2 x_1 - a_1 x_2 ) } \, ,
\end{equation}
for each $ {\xb} = (x_1, x_2) \in
 \left( \frac{1}{N} \mathbb Z \slash \mathbb Z \right)^2 $,
and one forms the series
\begin{equation} \label{basis}
  \phi_{\xb} (z) \, := \, (2 \pi N)^{-2} \, \sum_{{\ab} \in
 \left( \frac{1}{N} \mathbb Z \slash \mathbb Z \right)^2} \,
 \chi_{\xb} ({\ab}) \cdot  G_{\ab} (z) \, .
\end{equation}
 The definition is independent of $N$ and, for each $ {\xb} = (x_1, x_2) \in
 \left( \frac{1}{N} \mathbb Z \slash \mathbb Z \right)^2 \setminus
 {\0b} $ then $ \, \phi_{\xb} $, gives a weight $2$ Eisenstein series
 of level $N \,$.

\noindent To account for the special case when $ \, {\xb} =  {\0b}
\,$, one adjoins the non-holomorphic but modular function $
\, \phi_{\0b} \,$ defined in (\ref{gb0}).
\medskip

\noindent All the linear relations among the functions $\,
\phi_{\xb} \, , \, {\xb} \in ( \mathbb Q \slash \mathbb Z )^2  \,$
are encoded in the {\it distribution property}
\begin{equation} \label{dis}
  \phi_{\xb} \, = \, \sum_{{\yb} \cdot\check{\g} = {\xb}} \,
  \phi_{\yb} | \g  \, .
\end{equation}
where
$$
        \check{\g} \, = \, \det \g \cdot {\g}^{-1} \, .
$$
This allows to equip the extended Eisenstein space
$$
\Ec^*_2 (\mathbb Q) \, = \, \Ec_2 (\mathbb Q) \, \oplus \, \mathbb
Q \cdot \phi_0 \, ,
$$
with a linear $\, \text{PGL}^+ (2, \mathbb{Q})$-action, as
follows. Denoting
$$  \Sc \, := \, ( \mathbb Q \slash \mathbb Z )^2 \, ,
\quad \text{resp.} \quad  \, {\Sc}' \, := \, \Sc \setminus {\0b}
\,
$$
and identifying in the obvious way
$$ \text{PGL}^+ (2, \mathbb{Q}) \, \cong
\, M_2^+ (\mathbb{Z}) \slash \{ \text{scalars} \} \, ,
$$
where $\, M_2^+ (\mathbb{Z})  \, $ stands for the set of integral
$2 \times 2$-matrices of determinant $ > 0$, one defines the
action of $\, \g \in M_2^+ (\mathbb{Z}) \,$ by:
\begin{equation*} \label{a1}
  {\xb} | \g \, := \, \sum_{{\yb} \cdot\check{\g} = {\xb}} \, {\yb}
   \, \in \, \mathbb{Q} [\Sc] \, .
\end{equation*}
With this definition one has
\begin{equation*} \label{a2}
   \phi_{\xb} | \g \, = \, \phi_{{\xb} | \g} \, ,
   \qquad  \g \in M_2^+ (\mathbb{Z}) \, .
\end{equation*}
Modulo the subspace of  `distribution relations'
$$
 \Rc \, := \, \mathbb{Q}-\text{span of} \quad
 \left\{ {\xb} \, - \, {\xb} | \begin{pmatrix} n &0 \\
 0 &n  \end{pmatrix} \, \, ; \, {\xb} \in {\Sc} \, ,
 n \in \mathbb{Z} \setminus 0 \, \right\} \, ,
$$
the assignment $\, {\xb} \in {\Sc} \longmapsto \phi_{\xb} \,$
induces an isomorphism of $\, \text{PGL}_2^+ (2,
\mathbb{Q})$-modules
\begin{equation*} \label{miso}
 \mathbb{Q} [\Sc] / \Rc \, \cong \, \Ec^*_2 (\mathbb Q) \, .
\end{equation*}
\medskip

In view of the above the identity (\ref{gb0})  can be completed as follows:
\begin{equation}  \label{gb}
 \mu_{\g} \, = \,  2 ( \phi_{\0b} | \g \, - \, \phi_{\0b} ) \, =
 \, 2 \, \left( \sum_{{\yb} \cdot\check{\g} = {\0b}} \, \phi_{\yb}
 \, - \, \phi_{\0b} \right) \, , \quad \fl \,  \g \in M_2^+ (\mathbb{Z}) .
\end{equation}
\bigskip

\section{Characteristic map for the standard action}

In this section we provide the conceptual explanation for the period
pairing which was employed in \cite{CMmha} to obtain the Euler class
out of the universal Godbillon-Vey class $[\d_1] \in HC_{(\d, 1)}^1 \, (\Hc_1)$,
by showing that it is in fact the by-product of a characteristic map associated to
a $1$-trace which is invariant with respect to the standard action of
$\Hc_1$ on $\Ac$.

Extending a classical `splitting formula' for the restriction to
$\SL(2, \Zb)$ of the $2$-cocycle that gives the
universal cover of $\SL(2, \Rb)$,
we shall then obtain a new formula for the rational $2$-cocycle
representing the Euler class found in \cite{CMmha}, showing
that it differs from the Petersson cocycle  (\cite{As}, \cite{Pet}) by precisely
the coboundary of the classical Rademacher function.

\subsection{Characteristic map and cup products}

In~\cite{CMhti} (see also \cite{CMcchs}) we defined a
characteristic map associated to a Hopf module algebra with
invariant trace. The construction has been subsequently
extended to higher traces (cf.~\cite{Go})
and turned into a cup product
in Hopf-cyclic cohomology (cf.~\cite{KR}).
A predecessor of these constructions is the
contraction of a cyclic $n$-cocycle by the generator of a
$1$-parameter group of automorphisms that fixes the cocycle,
cf.~\cite[$Chap. III. \, 6. \b$]{Cncg}.
We shall apply the latter to a specific $1$-trace
$\tau_0 \in ZC^1 (\Ac)$, which will be described in
details in the next subsection. Further on, it
will also be applied in the context of
cyclic cohomology with coefficients, to $1$-traces
$\t_W \in ZC^1 (\Ac , W) $, where $W$ denotes an algebraically
irreducible $\GL(2, \Qb)$-module.
\medskip

Let us assume that $\t \in C^1(\Ac)$ is a cyclic cocycle which
 satisfies, with respect to a given Hopf action
 of $\Hc_1$ on $\Ac$, the invariance property
 \begin{equation} \label{invp0}
 \t (h_{(1)}(a^0) , h_{(2)}(a^1)) \, =
 \, \d(h) \, \t (a^0 , a^1) , \qquad \fl \, h \in \Hc_1 , \,  \,
 a^0 , a^1 \in \Ac \, .
 \end{equation}
The simplest expression for the cup product
$$
gv \, = \, \d_1 \# \t  \in ZC^2 (\Ac)
$$
is given by the contraction formula in~\cite[$Chap. III. \, 6. \b$]{Cncg}
mentioned above,
which (in the non-normalized form, cf.~\cite[\S 3]{Go}) takes
the expression:
\begin{equation} \label{gv}
gv \, (a^0 , a^1 , a^2) \, = \, \t (a^0 \, \d_1 (a^1), a^2) \, ,
\qquad a^0 , a^1 , a^2 \in \Ac \, .
 \end{equation}
 For the convenience of the reader, let
us check directly that this
formula gives a cocycle in the $(b, B)$-bicomplex of the
algebra $\Ac$.
\bigskip

\begin{lemma} \label{lgv}
Let $\t \in ZC^1(\Ac)$ be a cyclic cocycle
 satisfying the $\Hc_1$-invariance property (\ref{invp0}).
Then
$\quad b (gv )\, = \, 0 \quad \text{and} \quad B (gv) \, = \, 0 \,$.
\end{lemma}
\smallskip

\begin{proof} Using the fact that $\d_1$ acts as a derivation, one has
\begin{eqnarray} \nonumber
&\ &b (gv )\,(a^0 , a^1 , a^2 , a^3) \, =\, \t \, (a^0 a^1 \d_1 (a^2) , a^3 ) \, - \,
\t \, (a^0  \d_1 (a^1 a^2) , a^3 ) \cr \cr
&+& \t \, (a^0 \d_1 (a^1),  a^2 a^3 ) \, - \,
 \t \, (a^3 a^0 \d_1(a^1) , a^2 )  \cr \cr
&=& - \t \, (a^0 \d_1 (a^1)  a^2 , a^3 ) \, + \,
  \t \, (a^0  \d_1 (a^1),  a^2 a^3 ) \cr \cr
 &-& \t \, (a^3 a^0 \d_1(a^1) , a^2 )
\, = \,- b\t \, (a^0 \d_1 (a^1),  a^2 , a^3 ) \, = \, 0 \, .
\end{eqnarray}
Passing to $B$, since $\, \t (a, 1) \, = \, 0$ for any $a\in \Ac$, one has
\begin{eqnarray} \nonumber
&\,&B (gv )\,(a^0 , a^1) \, = \, gv \, (1, a^0, a^1) \, - \,
gv \, (1, a^1, a^0) \cr \cr
 &=& \t \, ( \d_1 (a^0) , a^1 ) \, - \,
\t \, (\d_1 (a^1) , a^0 )  \, = \,
\t \, ( \d_1 (a^0) , a^1 ) \, + \,  \t \, ( a^0 , \d_1 (a^1) ) \, = \, 0 \, ,
 \end{eqnarray}
the vanishing taking place because $\t \in ZC^1 (\Ac)$ is $\Hc_1$-invariant.
\end{proof}
\bigskip

\subsection{The basic invariant $1$-cocycles} \label{cusp}

The modular symbol cocycle
of weight $2$ associated to
a base point $z_0 \in \Hb$,  $ \t_0 \in ZC^1(\Ac)$,  is defined as follows.
For  monomials $f^0 U_{\g_0} , f^1 U_{\g_1}  \in \Ac$,
$\t_0 (f^0 U_{\g_0} , f^1 U_{\g_1} ) \, = \, 0$, unless they satisfy
the condition
 \begin{equation}  \label{2}
  w(f^0) + w(f^1) = 2 \qquad \text{and} \qquad
 \g_0  \g_1 =1 \, ,
 \end{equation}
in which case it is given by the integral
\begin{equation} \label{tau0}
\t_0 (f^0 U_{\g_0} , f^1 U_{\g_1} ) \, = \, \int_{z_0}^{ \g_0 \, z_0}
 f^0 \, f^1 | {\g_1}  \, dz \, .
 \end{equation}
 The fact that $\t_0 \in C^1(\Ac)$ is indeed an $\Hc_1$-invariant
 cyclic cocycle is the content of the following result.
\medskip

 \begin{proposition} \label{invt0}
The cochain $\t_0 \in C^1(\Ac)$ is a cyclic cocycle which
 satisfies the $\Hc_1$-invariance property (\ref{invp0}) with respect
 to the standard action.
 \end{proposition}
  \smallskip

 \begin{proof}
 Let $f^0 U_{\g_0} , f^1 U_{\g_1} , f^2 U_{\g_2}  \in \Ac$ be such that
  \begin{equation}  \label{3}
  w(f^0) + w(f^1)  + w(f^2) = 2 \qquad \text{and} \qquad
 \g_0  \g_1   \g_2 =1 \
 \end{equation}
 One has
 \begin{eqnarray} \nonumber
&\ &b \t_0 (f^0 U_{\g_0} , f^1 U_{\g_1} , f^2 U_{\g_2} )  \, = \cr \cr \nonumber
&=&\t_0 (f^0 f^1| {\g_0}^{-1}  U_{\g_0 \g_1} , f^2 U_{\g_2} ) -
  \t_0 ( f^0 U_{\g_0} , f^1 f^2 | {\g_1}^{-1}  U_{\g_1 \g_2} )  \cr \cr  \nonumber
&+& \t_0 (f^2 f^0 | {\g_2}^{-1}  U_{\g_2 \g_0} , f^1 U_{\g_1} )  \cr \cr \nonumber
&=&\int_{z_0}^{\g_0 \g_1 \, z_0}  f^0 f^1| {\g_0}^{-1}  f^2 | {\g_1}^{-1} {\g_0}^{-1} dz
-  \int_{z_0}^{\g_0  \, z_0}  f^0 f^1| {\g_0}^{-1}  f^2 | {\g_1}^{-1} {\g_0}^{-1} dz
 \cr \cr \nonumber
 &+&\int_{z_0}^{\g_2 \g_0 \, z_0}  f^2 f^0 | {\g_2}^{-1}  f^1 | {\g_0}^{-1} {\g_2}^{-1} dz
 \cr \cr \nonumber
 &=&\int_{\g_0 z_0}^{\g_0 \g_1 z_0}  f^0 f^1| {\g_0}^{-1}  f^2 | {\g_2} dz
 +  \int_{z_0}^{\g_2 \g_0 \, z_0}
  f^2 f^0 | {\g_2}^{-1}  f^1 | {\g_0}^{-1} {\g_2}^{-1} dz  \cr \cr \nonumber
  &=&\int_{\g_0 \, z_0}^{{\g_2}^{-1} \, z_0}  f^0 f^1| {\g_0}^{-1}  f^2 | {\g_2} dz
  +  \int_{{\g_2}^{-1} \, z_0}^{ \g_0 \, z_0}
    f^2 | \g_2  f^0  f^1 | {\g_0}^{-1} dz  \quad = \quad 0 \, ,
   \end{eqnarray}
 and so $ \t_0 $ is a Hochschild cocycle.

It is also cyclic, because for $f^0 U_{\g_0} , f^1 U_{\g_1}  \in \Ac$ satisfying
(\ref{2}) one has
  \begin{eqnarray} \nonumber
 \lb_1 \t_0 (f^0 U_{\g_0} , f^1 U_{\g_1} ) &=&
  -  \t_0 (f^1 U_{\g_1} , f^0 U_{\g_0} ) =
    -  \int_{z_0}^{ \g_1 \, z_0}  f^1 f^0 | {\g_1}^{-1}  dz   \cr \cr \nonumber
  &=&  \int_{z_0}^{ \g_0 \, z_0}  f^1 | {\g_1} f^0   dz
\,   = \, \t_0 (f^0 U_{\g_0} , f^1 U_{\g_1} ) \, .
  \end{eqnarray}

In view of its multiplicative nature, it suffices to check the $\Hc_1$-invariance
property (\ref{invp0}) on the algebra generators $\{Y, X, \d_1 \}$.
Starting with $Y$,  and with
$f^0 U_{\g_0} , f^1 U_{\g_1} $ satisfying (\ref{2}),
one has
 \begin{eqnarray} \nonumber
&\ &  \t_0 (Y(f^0 U_{\g_0}) , f^1 U_{\g_1} ) +
 \t_0 (f^0 U_{\g_0} , Y(f^1 U_{\g_1} )) =
\cr \cr \nonumber
&=&  \int_{z_0}^{ \g_0 \, z_0} Y( f^0 f^1 | {\g_0}^{-1})  dz
= \frac{w(f^0) + w(f^1)}{2} \int_{z_0}^{ \g_0 \, z_0}  f^0 f^1 | {\g_0}^{-1}  dz
\cr \cr \nonumber
&=& \int_{z_0}^{ \g_0 \, z_0}   f^0 f^1 | {\g_0}^{-1}  dz \quad = \quad
\d(Y) \, \t_0 (f^0 U_{\g_0} , f^1U_{\g_1} ) \, .
\end{eqnarray}
 Passing to $X$, the identity  (\ref{invp0}) is nontrivial only if
   $f^0 U_{\g_0} , f^1 U_{\g_1}  \in \Ac$ satisfy
 \begin{equation}  \label{00}
 w(f^0) + w(f^1) = 0 \qquad \text{and} \qquad
 \g_0  \g_1 =1 \, ,
 \end{equation}
 which actually implies that $f^0$ and $f^1$ are constants. One gets
\begin{equation}\nonumber
 \t (X(f^0) U_{\g_0} , f^1 U_{\g_1} ) +  \t (f^0 U_{\g_0} ,
X(f^1) U_{\g_1} ) +  \t (\d_1(f^0 U_{\g_0}) , Y(f^1) U_{\g_1} ) =
0
\end{equation}
since $X(f^j)=0$ and $Y(f^1)=0$.
 \smallskip

 Finally, with $f^0 U_{\g_0} , f^1 U_{\g_1} $ as above, one has
   \begin{eqnarray} \nonumber
&\ &  \t (\d_1(f^0 U_{\g_0}) , f^1 U_{\g_1} ) +
 \t (f^0 U_{\g_0} , \d_1(f^1 U_{\g_1} )) =
\cr \cr \nonumber
&=& \int_{z_0}^{ \g_0 \, z_0} \left( \mu_{\g_1}  +
  \mu_{\g_0}  |  {\g_1} \right)   f^0 f^1 | {\g_1}   dz \, = \, 0 \, ,
 \end{eqnarray}
because of the cocycle property of $\mu$.
 \end{proof}
\bigskip

Allowing the base point $z_0$ to belong to the `arithmetic' boundary
of the upper half plane $P^1(\Qb)$ requires some regularization of the
integral. This can be achieved by the standard procedure of removing
the poles of Eisenstein series (cf. \cite{St}). In the case at hand, it
amounts to a coboundary modification which we proceed now to describe.
\smallskip

 To obtain it, we start from the observation that
 the derivative of $\t_0$ with respect to the base point $z_0 \in \Hb$
  is a coboundary:
 \begin{equation}   \nonumber
\frac{d}{d z_0} \t_0 (f^0 U_{\g_0} , f^1 U_{\g_1} ) =
 (f^0 | {\g_0} \,  f^1)(z_0) - (f^0 \, f^1 | {\g_1} ) (z_0)
= - b \e (f^0 U_{\g_0} , f^1 U_{\g_1} ) ,
 \end{equation}
where $\e$ is the evaluation map at $z_0$:
 \begin{equation}   \nonumber
\e (f \, U_{\g} ) =  \left\{ \begin{matrix} f (z_0) \qquad \text{if}
\quad \g = 1 \,,\; f \in \Mc_2 \cr \cr 0 \qquad \qquad
\text{otherwise} \, .
\end{matrix} \right.
 \end{equation}

Taking the base point at the cusp $\infty$, we split the evaluation
functional $\e$ into two other functionals,
the constant term at $\infty$
 \begin{equation}   \nonumber
{\ab}_0  (f \, U_{\g} ) =  \left\{ \begin{matrix} a_0 \qquad
\text{if} \quad \g = 1 \,,\; f \in \Mc_2  \cr \cr 0 \qquad \qquad
\text{otherwise}
\end{matrix} \right.
 \end{equation}
and the evaluation of the remainder
  \begin{equation}   \nonumber
\wt{\e} (f \, U_{\g} ) =  \left\{ \begin{matrix}
 \wt{f} (z_0) \qquad \text{if} \quad \g = 1 \,,\; f \in \Mc_2 \cr \cr
 0 \qquad \qquad \text{otherwise} \, ,
\end{matrix} \right.
 \end{equation}
where  for $ \, f \in \Mc_2 \, $ of level $N$,
$$
f(z) \, = \, \sum_{n=0}^{\ify} \, a_n  \,e^{\frac{2\pi i n z}{N}}
$$
represents its Fourier expansion at $\ify$, and
$$
 \wt{f} (z) \, := \, f(z) \, - \, \ {\ab}_0 (f) \, \, .
$$
Both functionals are well-defined, because $ a_0 \,$  is
independent of the level.
To obtain a cohomologous cocycle
independent of $z_0$, it suffices to add to
$\t_0$  the sum of coboundaries of suitable
anti-derivatives for the two components. We therefore define
\begin{equation} \nonumber
\wt{\t_0} \,  = \,  \t_0  \, +  \, z_0  \, b\,  {\ab}_0 \, - \,
\int_{z_0}^{i \ify}  b\,  \wt{\e} \, dz \, ,
 \end{equation}
that is, for $f^0 U_{\g_0} , f^1 U_{\g_1}  \in \Ac$
as in (\ref{2})
\begin{eqnarray} \label{tt2}
\wt{\t_0} (f^0 U_{\g_0} , f^1 U_{\g_1} ) &=&  \int_{z_0}^{ \g_0 \, z_0}
 f^0\,  f^1 | {\g_1} \, dz \, + \, z_0 \, {\ab}_0 ( f^0 \, f^1 | {\g_1}  - f^0 | {\g_0} \,  f^1) \cr \cr
 &-& \int_{z_0}^{i \ify} \wt{\e} (f^0 \, f^1 | {\g_1}  - f^0 | {\g_0} \,  f^1) \, dz  \, .
 \end{eqnarray}

 The fact that $\, \wt{\t_0}  \in C^1(\Ac) \,$ still satisfies
Proposition \ref{invt0} can be checked by an obvious adaptation of its proof.
\bigskip

\subsection{Euler class and the transgression formula}

As a first example, we now specialize the construction of the cup product
to the invariant cocycle $\t_0 \in ZC^1 (\Ac)$.
Let $f^0 U_{\g_0} , f^1 U_{\g_1} , f^2 U_{\g_2}  \in \Ac$ be such that
$\,  \g_0  \g_1  \g_2 =1$. Then
 \begin{eqnarray} \label{Gv}
&\ & gv \,(f^0 U_{\g_0} , f^1 U_{\g_1} , f^2 U_{\g_2} ) \, = \,
\t_0 \, (f^0 U_{\g_0} \mu_{{\g_1}^{-1}} f^1 U_{\g_1} , f^2 U_{\g_2} ) =
  \cr \cr
&=& \t_0 (f^0 f^1| {\g_0}^{-1} \mu_{{\g_1}^{-1}} | {\g_0}^{-1} U_{\g_0 \g_1} ,
f^2 U_{\g_2} ) \cr \cr
&=& \int_{z_0}^{\g_0 \g_1z_0} f^0 f^1| {\g_0}^{-1} \mu_{{\g_1}^{-1}} | {\g_0}^{-1}
f^2 | {\g_1}^{-1} {\g_0}^{-1} \, dz \cr \cr
&=&  \int_{ {\g_0}^{-1} z_0}^{ \g_1 z_0} f^0 | \g_0 \, f^1 \,  \mu_{{\g_1}^{-1}}
\, f^2 | {\g_1}^{-1}  \, dz
=  \int_{ \g_2 z_0}^{ z_0} f^0 | \g_0 \g_1 \, f^1 | \g_1 \,
 \mu_{{\g_1}^{-1}} | \g_1 \, f^2  \, dz \cr \cr
 &=&  \int_{ z_0}^{\g_2  z_0} f^0 | \g_0 \g_1 \, f^1 |\g_1  \,
 f^2 \,  \mu_{\g_1}  \, dz \, .
\end{eqnarray}
In particular, the restriction to $\Ac_0 = \Cb [\GL^+ (2, \Qb)]$ is the group
cocycle
\begin{equation} \label{GV}
GV\, (\g_1, \g_2) := \, gv \, (U_{\g_0} , U_{\g_1} , U_{\g_2} ) \, = \,
\int_{ z_0}^{\g_2  z_0}   \mu_{\g_1}  \, dz \, ,
\end{equation}
whose real part
\begin{equation} \label{RGV}
\Re GV\, (\g_1, \g_2) := \, \Re \int_{ z_0}^{\g_2  z_0}   \mu_{\g_1}  \, dz  \, , \qquad
 \g_1 , \g_2 \in  \GL^+ (2, \mathbb{Q})
\end{equation}
represents a generator of $H^2 (\SL (2, \Qb), \Rb)$, hence a multiple
of the Euler class (cf. \cite[Thm. 16]{CMmha}).

For a more precise identification, we shall be very specific about
the choice of the Euler class. Namely, we take it as the class $\, {\eb}
\in H_{\rm bor}^2 (\Tb, \Zb)$ defined by the extension
\begin{equation} \label{eu}  \nonumber
0 \ra \Zb \ra \Rb \ra \Tb \ra 1 \, ;
\end{equation}
via the canonical isomorphisms
$$
H_{\rm bor}^2 (\Tb, \Zb) \simeq H^2 (B\Tb, \Zb)  \simeq H^2 (B\SL(2, \Rb), \Zb)
\simeq H_{\rm bor}^2 (\SL(2, \Rb), \Zb)
$$
followed by the succession of natural map
$$
H_{\rm bor}^2 (\SL(2, \Rb), \Zb)
\ra H^2 (\SL (2, \Qb), \Zb) \ra H^2 (\SL (2, \Qb), \Rb)
$$
we regard it as a class $\, {\eb} \in H^2 (\SL (2, \Qb), \Rb) $.
\bigskip

\begin{proposition} \label{ptr}
    The $2$-cocycle $ \Re GV
\in \, Z^{2} (\SL (2, \mathbb{Q}), \mathbb{R}) $
represents the class $\,- 2{\eb} \, \in \, H^{2}(\SL (2, \Qb), \Rb)$,
while $\Im GV$ is a coboundary.
\end{proposition}
\medskip

\begin{proof} In view of the definition (\ref{GV}) and using
(\ref{mu}), one has for
$ \g_1 , \g_2 \in  \GL^+ (2, \mathbb{Q})$,
\begin{eqnarray} \label{1t}
12  \pi i \, GV (\g_1, \g_2) &=&12  \pi i \, \int_{z_0}^{\g_{2}  z_0} \,
\mu_{\g_1} (z) \, dz \, = \, \int_{z_0}^{\g_{2}  z_0} \,
\frac{d}{dz} \log \frac{\D | \g_1}{\D} \, dz \cr \cr
 &=& \, \int_{z_0}^{\g_{2}  z_0} \,
(d \log \D | \g_1- \,d \log \D)   \cr \cr &=& \log \D | \g_1 (
\g_{2} z_0 ) - \log \D | \g_1 ( z_0 )  \cr \cr &-& \left( \log \D (
\g_{2} z_0 ) - \log \D ( z_0 ) \right)
\end{eqnarray}
where, since both $\D $ and $\D | \g_1$ don't have zeros in $\mathbb
H$ one lets  $\log \D $ and $\log \D | \g_1$ be holomorphic
determinations of the logarithm whose choice is unimportant at this
stage, since the additive constant which depends only on $\g_{1}$
cancels out. Let
$$
j(\g, z) \, = \, cz + d \, , \qquad  \g = \begin{pmatrix} a &b\\
 c &d  \end{pmatrix} \in \GL^+ (2, \mathbb{Q}) \, ,
$$
be the automorphy factor. Since it has no zero in $\mathbb H$ one
can choose for each $\g$ a holomorphic determination $\log j^{2}
(\g, z )$ of its logarithm (for instance using the principal branch
for $\, \mathbb{C} \setminus [0, \infty) \,$ of the logarithm when
$c\neq 0$ and taking $\log d^{2}$ when $c=0$). One then has,
$$
\log \D | \g ( z ) \, = \, \log \D (\g  z ) \, - \, 6 \, \ \Log
j^{2} (\g,  z ) \, + \, 2 \pi i \, k (\g) \, ,\quad \forall z \in
\mathbb H
$$
for some $\,  k (\g) \in \mathbb{Z} \,$. Thus (\ref{1t}) can be
continued as follows:
\begin{eqnarray}   \label{2t}
12 \pi i GV (\g_1, \g_2) &=&
 \log \D (\g_ 1 \g_2 z_0 ) - \log \D (\g_1 z_0 )
 - \log \D (\g_2  z_0 )  \cr \cr
 &+& \log \D ( z_0 )
-  6 \left( \log j^{2} (\g_1,  \g_2  z_0 ) - \log j^{2}
(\g_1, z_0 ) \right) \qquad \qquad
\end{eqnarray}
again after the cancelation of the additive constants. The equality
\begin{equation} \label{log}
 \log j^{2} (\g_1 \g_2 ,  z_0 ) =  \log j^{2} (\g_1,  \g_2  z_0 )
 +  \log j^{2} (\g_2 ,  z_0 ) - 2 \pi i \, c (\g_1,  \g_2 ) \, ,
\end{equation}
determines a cocycle $c \in Z^2 (\PSL (2, \mathbb{R}), \,
\mathbb{Z}) \,$  (which is precisely the cocycle discussed
in~\cite[\S B-2]{BG}, and whose cohomology class is independent of
the choices of the branches $\log j^{2} (\g, z )$ of the logarithm).
Inserting (\ref{log}) into (\ref{2t}) one obtains
\begin{eqnarray} \label{3t}
12 \pi i GV (\g_1, \g_2) &=&  \log \D ( z_0 )
+ \log \D (\g_ 1 \g_2  z_0 ) - \log \D (\g_1 z_0 )
 - \log \D (\g_2  z_0 ) \cr \cr
&-& 6 \left( \log j^{2} (\g_1 \g_2 ,  z_0 ) - \log j^{2} (\g_2,
z_0 )
- \log j^{2} (\g_1, z_0 ) \right) \cr \cr
&-& 12 \pi i \, c (\g_1,  \g_2 ) \, .
\end{eqnarray}
This identity shows that the cocycles $ \Re GV $ and $\, -c $ are
cohomologous in $Z^2 ( \PSL(2, \Qb), \Rb)$, and also
that $ \Im GV $ is a coboundary.

On the other hand, the restriction of
 $c \in Z_{\rm bor}^2 (\SL(2, \Rb), \Zb)  $ to $\Tb = SO(2)$,
\begin{eqnarray} \nonumber
c(\g(\th_1),  \g(\th_2)) &=& \frac{1}{2 \pi i} \left( \Log e^{2i
\th_1} + \Log e^{2i \th_2} - \Log e^{2i (\th_1 + \th_2)} \right) \,
,  \cr \cr \quad \text{where} \quad
\g(\th) &=&  \begin{pmatrix} \cos \th &-\sin \th\\
 \sin \th &\cos \th  \end{pmatrix} \, , \qquad \th \in [0, 2 \pi) \, .
\end{eqnarray}
evidently represents the class $2 {\eb} \in H_{\rm bor}^2 (\Tb, \Zb)$,
which concludes the proof.
 \end{proof}
 \bigskip

 The Euler class ${\eb} \in H^2 (\SL(2, \Qb), \Zb)$
 occurs naturally in the context of
 the Chern character in K-homology \cite{C1},
 as the Chern character of a natural Fredholm module
 given by the "dual Dirac" operator relative to a
 base point $z_0 \in \Hb$. When moving the base
 point to the cusp $ i \infty \in P^1(\Qb)$, it coincides with the
  restriction of the $2$-cocycle $ e \in Z^2 (\SL(2, \Rb), \Zb)$
introduced by Petersson (cf. \cite{Pet}) and
investigated in detail by Asai (cf.  \cite{As}, where it is denoted $w$).
It is defined, for $g_1, g_2 \in \SL(2, \Rb )$, by the formula
\begin{equation} \label{as}
e( g_1, g_2) \, = \, \frac{1}{2 \pi i}
\left(\log  j (g_2 , z) +  \log  j (g_1 , g_2 z) - \log  j (g_1 g_2 , z)\right) \, ,
\end{equation}
with the logarithm chosen so that $\Im \log \in [- \pi, \pi)$; the
above definition is independent of $z \in \Hb$.

Asai  \cite[\S 1-4]{As} has shown that it can be given
a simple expression,
analogous to Kubota's cocycles \cite{Ku} for coverings over local fields,
which is as follows:
\begin{equation} \label{hs}
e( g_1, g_2) \, = \, - (x(g_1)| x(g_2) ) \, + \, (- x(g_1) x(g_2) | x(g_1 g_2)) \, ,
\end{equation}
where
\begin{equation} \nonumber
 x(g) \, = \left\{\begin{matrix} &c& \quad \text{if} \quad c>0 \, , \cr \cr
 &d& \quad  \text{if} \quad c=0 \, , \end{matrix} \right.
 \qquad \fl \, g =  \begin{pmatrix} a &b\\
 c &d  \end{pmatrix} \in  \SL(2, \Rb )
 \end{equation}
and for any two numbers $x_1 , x_2 \in \Rb$, the (Hilbert-like) symbol
$(x_1 | x_2)$ is defined as
\begin{equation} \nonumber
(x_1 | x_2) \, =  \left\{\begin{matrix} &1& \qquad \text{iff} \quad
x_1 <0 \quad \text{and} \quad x_2 <0  \, ,\cr \cr
 &0& \qquad  \text{otherwise} \, . \end{matrix} \right.
 \end{equation}
\bigskip

Replacing $\t_0$ by $\wt{\t}_0$ one obtains, as in (\ref{Gv}) and (\ref{GV}),
the cohomologous group $2$-cocycle on $\GL^+ (2, \mathbb{Q}) $
\begin{eqnarray} \label{tGV}
&\ &\wt{GV}\, (\g_1, \g_2) \,= \, \wt{\t}_0 (U_{\g_0} \, \d_1(U_{\g_1}) , U_{\g_2}) \, =
\, \wt{\t}_0 (\mu_{{\g_1}^{-1}} | \g_0^{-1} \, U_{\g_0 \g_1} , U_{\g_2} )
\qquad \qquad \cr \cr
&=& \int_{ z_0}^{\g_2  z_0}   \mu_{\g_1} dz
\, +\, z_0\, {\ab}_0 ( \mu_{\g_1} - \mu_{\g_1} | \g_2 )
+ \int_{ z_0}^{i \infty} (\wt{\mu_{\g_1} | \g_2} - \wt{\mu_{\g_1}}) dz  \, .
\end{eqnarray}
In \cite[\S 4]{CMmha} we found an explicit rational formula for
$\Re \wt{GV}$,
in terms of Rademacher-Dedekind sums, which we proceed now to recall.
\medskip

Since by its very definition $\Re \wt{GV}$ descends to a
$2$-cocycle on ${\rm PGL}^+ (2, \mathbb{Q}) $, it suffices to express it
for pairs of matrices with integer entries
\begin{equation} \label{gg}
 \g_1 \, = \, \begin{pmatrix} a_1 &b_1 \\
 c_1 &d_1  \end{pmatrix} \, , \quad
\g_2 \, = \, \begin{pmatrix} a_2 &b_2\\
 c_2 &d_2  \end{pmatrix} \,  \in \, M_2^+ (\mathbb{Z})  \, .
\end{equation}

When  $ \g_2 \, \in \, B^+ (\mathbb{Z})$, that is
$\,  c_2 \, = \, 0$, then
\begin{equation} \label{f1}
  \Re \wt{GV} (\g_1, \g_2) \, = \,  \frac{b_2}{d_2} \,
  \sum_{{\xb} \cdot \check{\g_1} = {\0b}, {\xb} \neq {\0b}} \,  {\Bb}_{2} (x_1)
  \, ,
\end{equation}
where $\, {\Bb}_2 (x):= (x-[x] )^2 - (x-[x]) + \frac{1}{6} \,$,
 with $[x] =$  greatest
 integer $ \leq x $.
 \smallskip

When $c_2  > 0$, then
 \begin{eqnarray} \label{f2}
&\ &\Re \wt{GV} (\g_1, \g_2) \,  = \, \frac{a_2}{c_2} \,  \sum_{{\xb} \cdot
\check{\g_1} = {\0b} , {\xb} \neq {\0b}}  {\Bb}_{2} (x_1)
\, + \, \frac{d_2}{c_2}\,   \sum_{{\xb} \cdot
 \check{\g_2} \check{\g_1} = {\0b} , {\xb} \neq {\0b}}
{\Bb}_{2} (x_1) \cr \cr
&-& 2 \sum_{{\xb} \cdot \check{\g_1} = {\0b} , {\xb} \neq {\0b}} \sum_{j=0}^{c'_2 -1}
 {\Bb}_1 \left(\frac{x_1 + j}{c'_2} \right)
 {\Bb}_1 \left(\frac{a'_2 (x_1 + j)}{c'_2} + x_2 \right)
  \end{eqnarray}
 where $\,  \frac{a'_2}{c'_2} = \frac{a_2}{c_2} \,$ , $\, (a'_2 , c'_2) = 1 \,$, and
$ {\Bb}_1 (x):= x-[x]  - \frac{1}{2} \,$, for any $\, x \in \Rb$ .
\bigskip

We shall obtain below a simpler expression for $\Re \wt{GV} \in
Z^2 ( \SL(2, \Qb), \Qb)$, through a transgression formula
which involves only the classical Dedekind sums, through
the Rademacher function.
Let us recall that,
for a pair of integers $m, n$ with $(m, n) = 1$ and $n \geq 1$,
the Dedekind sum is given by the formula
\begin{equation} \label{Rs}
 s(\frac{m}{n}) \, = \, \sum_{j=1}^{n-1}  {\Bb}_1 \left(\frac{j}{n} \right) \,
 {\Bb}_1 \left(\frac{m \, j}{n} \right) \, .
 \end{equation}
The Rademacher function $\Phi: \SL(2, \Zb) \ra \Zb$
is uniquely characterized (cf. \cite{RG})
by the coboundary relation
\begin{equation} \label{TR}
 \Phi(\s_1 \s_2)  =  \Phi(\s_1) - {\Phi}(\s_2) - 3 \sign(c_1 c_2 c_3) \, , \quad
 \s_1 , \s_2 \in \SL(2, \Qb)
 \end{equation}
 where  $\s_3 = \s_1 \s_2$, $\, \, \s_i =  \begin{pmatrix} a_i &b_i\\
 c_i &d_i  \end{pmatrix} \in \SL(2, \Zb )$, $i = 1, 2, 3 , \,$
and is explicitly given  by  the following formula  (cf. \cite{Ra}):
\begin{equation} \label{Rad}
\Phi (\s) \, =  \left\{\begin{matrix}  &\frac{b}{d}&  \qquad \qquad
\qquad \qquad \qquad \text{if} \quad
c = 0 ,\cr \cr
 &\frac{a+d}{c}& - \quad 12 \sign(c) s(\frac{a}{|c|}) \qquad \text{if} \quad c \neq 0
 \end{matrix} \right.
 \end{equation}
   for any $\s =  \begin{pmatrix} a &b\\
 c &d  \end{pmatrix} \in  \SL(2, \Zb )$.

 We extend it, in a slightly modified version,
 to a function $\wt{\Phi}: \GL^+(2, \Qb) \ra \Qb$, as follows. First, for any
 $\s =  \begin{pmatrix} a &b\\
 c &d  \end{pmatrix} \in \SL(2, \Zb )\,$, we define
\begin{equation} \label{Rads}
\wt{\Phi} (\s) \, =  \left\{\begin{matrix}\frac{b}{12 d} &+& \frac{1-\sign(d)}{4}
\qquad  \qquad \qquad \quad   \text{if} \quad
c = 0 ,\cr \cr
 \frac{a+d}{12 c} &-&  \sign(c)
 \left(\frac{1}{4} + s(\frac{a}{|c|})\right) \qquad \text{if} \quad c \neq 0 ,
 \end{matrix} \right.
 \end{equation}
while for any $\b \in B(\Qb)$, where
$ B(\Qb) = \left\{ \b =  \begin{pmatrix} a &b\\
 0 &d  \end{pmatrix} \in \GL^+(2, \Qb )\right\} \,$, we set
 \begin{equation} \label{Radb}
 \wt{\Phi} (\b) \, = \, \frac{b}{12 d} \, + \, \frac{1-\sign(d)}{4}  .
 \end{equation}
Now given $\g \in \GL^+(2, \Qb)$, after factoring it in the form
$$ \g \, = \, \s \cdot \b \qquad \text{with} \qquad
\s \in \SL(2, \Zb ) \quad  \text{and}  \quad \b \in B(\Qb) \, ,
$$
we define
 \begin{equation} \label{Radq}
 \wt{\Phi} (\g) \, = \,   \wt{\Phi} (\s) \, +  \wt{\Phi} (\b) \, ;
 \end{equation}
one easily checks, by elementary calculations, that the definition is consistent.
 \medskip

The transgression formula within the Euler class can now
be stated as follows.
 \medskip

\begin{theorem} \label{Tprt}
The function $\wt{\Phi}: \SL(2, \Qb) \ra \Qb$ is uniquely characterized
by the identity
\begin{equation} \label{TFE}
\frac{1}{2} \Re \wt{GV} (\g_1 , \g_2) \, + \, e (\g_1 , \g_2)   \, = \,
 \wt{\Phi}(\g_1 \g_2) \, - \,  \wt{\Phi}(\g_1)\,  - \, \wt{\Phi}(\g_2) \, ,
 \end{equation}
 for any $\g_1 , \g_2 \in \SL(2, \Qb)$.
\end{theorem}
\medskip

\begin{proof} By Proposition \ref{ptr}, $\displaystyle \, [\frac{1}{2} \Re \wt{GV} +
e] = 0\, $  in $H^2 (\SL(2, \Qb) , \Rb)$. Since
$$H^1 (\SL(2, \Qb) , \Rb) \, = \, 0 \, ,
$$
there exists a unique function $\Psi: \SL(2, \Qb) \ra \Rb$ such that
\begin{equation} \label{ET}
\frac{1}{2} \Re \wt{GV} (\g_1 , \g_2) + e (\g_1 , \g_2)  =
 \Psi(\g_1 \g_2) -  \Psi(\g_1) -  \Psi(\g_2) \, .
 \end{equation}
Restricting to $ \SL(2, \Zb)$, and taking into account that
$\mu_\s = 0$ for all $\s \in \SL(2, \Zb)$, one obtains
\begin{equation} \nonumber
  e (\s_1 , \s_2)  =
 \Psi(\s_1 \s_2) -  \Psi(\s_1) -  \Psi(\s_2) \, , \qquad \s_1, \s_2 \in  \SL(2, \Zb) .
 \end{equation}
This is the splitting formula (\ref{TR}) which uniquely characterizes the restriction
(\ref{Rads}) of $ \wt{\Phi}$ to $ \SL(2, \Zb)$ (see \cite[Thm. 3]{As}), so that
\begin{equation} \label{ons}
\Psi (\s) = \wt{\Phi}(\s)\, , \qquad \fl \, \s \in \SL(2, \Zb)\, .
 \end{equation}
Furthermore, taking in (\ref{ET}) $\g_1 = \s \in \SL(2, \Zb)$ and
$\g_2 = \b \in B_1^+(\Qb) =  \left\{ \b =  \begin{pmatrix} a &b\\
 0 &d  \end{pmatrix} \in B( \Qb ) \, ; \, a > 0 \, , \, a\, d = 1 \right\} \,$,
one obtains
\begin{equation} \label{plus}
\Psi (\s \b) \,  = \, \wt{\Phi}(\s) + \Psi(\b) \, , \quad \fl \,  \s \in \SL(2, \Zb) \, , \,
\b \in B_1^+(\Qb) \, ,
 \end{equation}
 because $e(\s, \b) = 0$, by  \cite[Lemma 3]{As}.
 In particular, for $\s = -I$, one has
\begin{equation} \label{-}
\Psi (- \b) \,  = \, \wt{\Phi}(-I) + \Psi(\b) \, = \, \Psi(\b) + \frac{1}{2},
\qquad \fl \, \b \in B_1^+(\Qb) \, ,
 \end{equation}
 since $\displaystyle  \wt{\Phi}(-I) = \frac{1}{2}$, cf. \cite[Lemma 4]{As}.

 It remains to prove that $\Psi (\b) =  \wt{\Phi}(\b)$ for any $\b \in B_1^+(\Qb)$.
 Specializing (\ref{ET}) to $B_1^+(\Qb)$, and recalling that $e$ vanishes on
  $B_1^+(\Qb)$, one obtains
 \begin{equation} \nonumber
 \Psi(\b_1 \b_2) -  \Psi(\b_1) -  \Psi(\b_2) \, = \,
 \frac{1}{2} \Re \wt{GV} (\b_1 , \b_2)\,
 \end{equation}
 and by (\ref{tGV}) and (\ref{gb0}) this can be computed as the real part of
 \begin{eqnarray} \nonumber
 &\, & \int_{ z_0}^{\b_2  z_0}
 ( \phi_0 | {\b_1} - \phi_0 )dz
\, +\, z_0\, {\ab}_0 (  \phi_0 | {\b_1} - \phi_0 -
\phi_0 | {\b_1}  \b_2  +  \phi_0 | \b_2)  \cr \cr
&+& \int_{ z_0}^{i \infty} (\wt{\phi_0 | {\b_1}  \b_2}  - \wt{\phi_0 | \b_2}
-  \wt{\phi_0 | \b_1} +  \wt{\phi_0}) dz   \cr \cr
 &=& \int_{ z_0}^{\b_2  z_0} ( \wt{\phi_0 | \b_1} -  \wt{\phi_0}) dz +
(\b_2  z_0 -  z_0 )\, {\ab}_0 (\phi_0 | \b_1 - \phi_0 ) \cr \cr
&+&  z_0\, {\ab}_0 (  \phi_0 | {\b_1} - \phi_0 -
\phi_0 | {\b_1}  \b_2  +  \phi_0 | \b_2)  \cr \cr
&+& \int_{ z_0}^{i \infty} (\wt{\phi_0 | {\b_1}  \b_2}  - \wt{\phi_0 | \b_2}
-  \wt{\phi_0 | \b_1} +  \wt{\phi_0}) dz \cr \cr
&=& \int_{ z_0}^{\b_2  z_0} ( \wt{\phi_0 | \b_1} -  \wt{\phi_0}) dz
+ \int_{ z_0}^{i \infty} (\wt{\phi_0 | {\b_1}  \b_2}  - \wt{\phi_0 | \b_2}
-  \wt{\phi_0 | \b_1} +  \wt{\phi_0}) dz \cr \cr
&+& \b_2  z_0 \, {\ab}_0 (\phi_0 | \b_1 - \phi_0 ) -
 z_0\, {\ab}_0 (  (\phi_0 | \b_1 - \phi_0 ) | \b_2) \, .
\end{eqnarray}
Since obviously
\begin{equation} \label{a0b}
{\ab}_0 (f | \b) \, = \, \frac{a}{d} \,  {\ab}_0 (f)\, , \quad \text{for} \quad
\b =  \begin{pmatrix} a &b\\
 0 &d  \end{pmatrix} \in B( \Qb ) \, ,
 \end{equation}
the last line contributes
 \begin{equation} \nonumber
\left(\frac{a_2\, z_0 + b_2}{d_2} \, -\,
 \frac{a_2}{d_2} \, z_0\right)  \, {\ab}_0 (\phi_0 | \b_1 - \phi_0 ) \, = \,
 \frac{b_2}{d_2}  \, {\ab}_0 (\phi_0 | \b_1 - \phi_0 ) \,
\end{equation}
which can be further computed as
$$
=  \, \frac{1}{12} \, \frac{b_2}{d_2}  \, \left( \frac{a_1}{d_1} - 1\right) \, ,
$$
by using (\ref{f1}) or more directly
the Fourier expansion of $\phi_0$ (cf. e.g. \cite[Prop. 2.4.2]{St}).
On the other hand, when $z_0 \ra i \ify$ then $\b_2 \, z_0 \ra i \ify$ too, so that
both integrals converge to $0$.
We conclude that
\begin{equation}  \label{BET}
 \Psi(\b_1 \b_2) -  \Psi(\b_1) -  \Psi(\b_2) \, = \,
\frac{1}{2} \Re \wt{GV} (\b_1 , \b_2)\,  = \,
   \frac{1}{12} \, \frac{b_2}{d_2}  \, \left( \frac{a_1}{d_1} - 1\right) \, .
 \end{equation}
Using (\ref{Radb}) it is elementary to check that, for $\b_1, \b_2 \in B_1^+( \Qb ) $,
\begin{equation} \nonumber
\wt{ \Phi}(\b_1 \b_2) -  \wt{\Phi}(\b_1) -  \wt{\Phi}(\b_2) \, = \,
  \frac{1}{12} \, \frac{b_2}{d_2}  \, \left( \frac{a_1}{d_1} - 1\right) \, .
 \end{equation}
Thus, $\Psi | B_1^+(\Qb)$ and $\wt{ \Phi} | B_1^+(\Qb)$
 can only differ by a character of $B_1^+(\Qb) $.

To show that they coincide, it suffices to prove that  $\Psi$ vanishes on the torus
$$
T =  \left\{  \d =  \begin{pmatrix} a &0\\
 0 &d  \end{pmatrix} \, ;  \, \d \in B_1^+( \Qb ) \right\} \, .
$$
For any $\d \in T$, one has
$$ \s_0 \, \d \, = \, \d^{-1} \, \s_0 \, , \quad \text{where} \quad
\s_0 =  \begin{pmatrix} 0 &-1\\
                   1 &\, \, 0  \end{pmatrix} \, ,
$$
hence
$$ \Psi (\s_0 \, \d ) \, = \, \Psi (\d^{-1} \, \s_0) \, .
$$
>From (\ref{ET}) it follows that
$$
\Psi (\s_0) + \Psi (\d) +  \frac{1}{2} \Re \wt{GV} (\s_0 , \d)\, = \,
\Psi (\d^{-1}) +  \Psi (\s_0) +  \frac{1}{2} \Re \wt{GV} (\d^{-1} , \s_0) \, ,
$$
therefore, since $\, \Re \wt{GV} (\s_0 , \d) =0$,
$$\Psi (\d)\, - \, \Psi (\d^{-1}) \, = \, \frac{1}{2} \, \Re \wt{GV} (\d^{-1} , \s_0)
 \, .
$$
On the other hand, by (\ref{BET})
$$
\Psi (\d) \, + \, \Psi (\d^{-1}) \, = \, 0 \, .
$$
Hence
$$
\Psi (\d) \,= \, \frac{1}{4} \, \Re \wt{GV} (\d^{-1} , \s_0) \, ,
$$
and it remains to show that the right hand side vanishes. We shall
apply formula (\ref{f2}), which allows us to replace $\d^{-1}$, after
multiplication by a scalar, with a diagonal matrix with positive integer entries
$ \rho \, =  \begin{pmatrix} m &0 \\
 0 &n  \end{pmatrix}   \in M_2^+ (\mathbb{Z})$. In this case, it simply gives
  \begin{eqnarray} \nonumber
\frac{1}{2 }\, \Re \wt{GV} (\rho, \s_0) &=&
-  \sum_{{\xb} \cdot \check{\rho} = {\0b}}   {\Bb}_1 (x_1) \, {\Bb}_1 ( x_2)
\,+  \,  {\Bb}_1 (0)^2  \cr \cr
&=& -  \sum_{j=0}^{m-1} {\Bb}_1 (\frac{j}{m}) \, \sum_{k=0}^{n-1} {\Bb}_1 (\frac{k}{n})
+   {\Bb}_1 (0)^2  \, = \,  0 \, ,
  \end{eqnarray}
because of the distribution property of the first Bernoulli function.
\end{proof}
\bigskip

\begin{remark}
{ \rm  Note that while the cocycle $e$ admits a $K$-homological interpretation as
a Chern character, the cohomologous cocycle $\frac{1}{2} \Re \wt{GV} $
obviously should also have such an interpretation. This in turns would
allow to put the above transgression on the same $K$-homological footing
as in \cite{C2}.}
\end{remark}

\section{ Modular symbol cyclic cocycles of higher weight}

We begin in this section to extend the above results to the
case of higher weight modular symbols.
This will involve introducing
the cyclic cohomology with coefficients
$\, HC^* (\Ac, W) \,$ and describing
the higher weight analogues of the
basic invariant cocycles. They
will be used in the next section
to define the characteristic maps
corresponding to the
degenerate  actions of weight $m \geq 2$ of $\Hc_1$ on
$\Ac= \Ac_{G^+ (\mathbb{Q})}$.

\subsection{Cyclic cohomology with coefficients}

We need to introduce the \textit{cyclic cohomology with coefficients}
$\, HC^* (\Ac, W) $, where $ \Ac$ is
the crossed product algebra
$\Ac_{G^+ (\mathbb{Q})} $  and $\, W \,$ is a
$\,  \GL^+ (2, \mathbb{Q}) \,$-module. It is a special case of the
Hopf-cyclic cohomology with coefficients, cf. \cite[\S 3]{HKRS}, in which the `gauge'
 Hopf algebra is the group ring
$$ \Gc \, = \, \mathbb{C} [G^+ (\mathbb{Q})] , \qquad G^+ (\mathbb{Q}) := \GL^+ (2, \mathbb{Q}) ,
$$
equipped with its usual Hopf algebra structure;
$W$ is viewed as a left $ \Gc $-module and as a trivial $ \Gc $-comodule, and
$\Ac $ is regarded as a left $ \Gc $-comodule algebra with respect to the
intrinsic coaction
$$ a = f \, U_{\g}  \longmapsto a_{(-1)} \ot a_{(0)} := U_{\g} \ot  f \, U_{\g} \in \Gc \ot \Ac \, .
$$
Thus, by definition, $\, HC^* (\Ac, W) \,$ is the cohomology associated to the cyclic module
\begin{equation} \label{tcc}
C^* (\Ac, W)  \, := \, \Hom^{\Gc} (\Ac^{*+1}, W) \, ,
\end{equation}
whose cyclic structure is defined by the operators
\begin{eqnarray}
\partial_i \phi (a^0 \ot \cdots \ot a^{n})&=& \, \phi (a^0 \ot \cdots a^i a^{i+1}
\ot \cdots \ot a^n ) \, , \,  0 \leq i < n \, , \nonumber \\
\partial_n \phi (a^0 \ot \cdots \ot a^{n})&=& \, S(a^n_{(-1)}) \, \phi (a^n_{(0)}a^0 \ot \cdots \ot a^{n-1})
 \, , \qquad \qquad \nonumber \\
\s_i \phi (a^0 \ot \cdots \ot a^{n})&=& \, \phi (a^0 \ot \cdots a^i \ot a^{i+1} \ot \cdots \ot a^{n}) \, ,
\, 0 \leq i < n \, , \nonumber \\
    \tau_n \phi (a^0 \ot \cdots \ot a^{n})&=& \, S(a^n_{(-1)}) \, \phi (a^n_{(0)} \ot a^0 \ot \cdots \ot a^{n-1})
    \, . \nonumber
\end{eqnarray}
\bigskip

\noindent For each $m \in \Nb$ we denote by
$\,W_m  $ the simple  $\SL (2, \mathbb{C})$-module of dimension $m+1$,
realized as the the space
$$
W_m = \{ P(T_1, T_2) \in \mathbb{C} [T_1, T_2] \, ; \, P \, \text{is homogeneous of degree} \, m  \}
$$
and we let $\GL ^+(2, \mathbb{R})$ act on $W_m  $ by
\begin{equation*}
 (g \, \cdot P) \, (T_1, T_2) \, = \, \det (g)^{- \frac{m}{2}} \, P(aT_1 + c T_2, bT_1 + dT_2) , \quad
g \, = \, \begin{pmatrix} a &b \\
 c &d  \end{pmatrix}  \, .
\end{equation*}
Note that as a $\GL ^+(2, \mathbb{R})$-module, $W_m  $ is the complexification
of
$$
W_m (\Rb) =
\{ P(T_1, T_2) \in \mathbb{R} [T_1, T_2] \, ; \, P \, \text{is homogeneous of degree} \, m  \} \, ,
$$
and we denote by $\, \Re : W_m \ra W_m (\Rb) $ the projection obtained by taking
the real parts of the coefficients. We also note that
as a  $\GL ^+(2, \mathbb{Q})$-module $W_m (\Rb)$
has an obvious rational structure
$$
W_m (\Qb) =
\{ P(T_1, T_2) \in \mathbb{Q} [T_1, T_2] \, ; \, P \, \text{is homogeneous of degree} \, m  \} \, .
$$
\smallskip

We denote by $\, F_m: \mathbb{H} \ra W_m \,$ the polynomial
function
\begin{equation}
F_m (z) \, = \,  (z T_1 \, + \, T_2)^m \, ,
\end{equation}
and note that it
satisfies, for any $g \in GL ^+(2, \mathbb{R})$, the covariance property
\begin{equation} \label{Fcov}
 F_m  \vert g \, (z) \equiv \det (g)^{- \frac{m}{2}}  (cz + d)^m \, F_m (g z)
  = g \cdot F_m  (z) \, .
  \end{equation}
\bigskip

\subsection{Cyclic $1$-cocycles with coefficients: base point in $\Hb$}

\noindent With these ingredients at hand, and after making the
additional choice of a `base point' $\, z_0 \in \mathbb{H}\,$, we
proceed to define the \textit{invariant cocycles} which will support
characteristic maps associated to  Hopf actions in the degenerate
case.

Regarding  $\, \Hom_{\Gc} (\Ac \ot \Ac , W_m) \,$ as a graded linear space
 with respect to the weight filtration inherited from $\, \Ac \,$, we define the
weight $2$ element  $\, \t_m \in \Hom_{\Gc} (\Ac \ot \Ac  , W_m) \,$ as follows.
Let
$$
 a^0 = \sum_{\a} f^0_{\a} U_{\a}   \in \Ac_{w(a^0)}  \quad \text{and}
 \quad a^1=\sum_{\b} f^1_{\b} U_{\b}  \in \Ac_{w(a^1)}  \, ,
 $$
 be two homogeneous elements in $  \Ac $; by definition,
\begin{equation} \nonumber
\t_m (a^0 , a^1 )  =
\left\{ \begin{matrix} \displaystyle
  \sum_{\a} \int_{z_0}^{\a \, z_0}  F_m \, f^0_{\a} \,  f^1_{\a^{-1}} \vert  \a^{-1}  \, dz  \quad
  \text{if}  \quad w(a^0 ) + w(a^1) = m+2, \cr \cr
  0 \qquad \qquad \qquad \qquad  \text{otherwise} \, .
\end{matrix} \right.
\end{equation}
  \bigskip

  \begin{lemma} \label{coco}
For each $m \geq 2$,
$ \t_{m-2} \in C^1 (\Ac , W_{m-2}) \,$ is a cyclic cocycle.
 \end{lemma}
 \smallskip

 \begin{proof} For notational convenience, we shall omit the subscript $m-2$ in
 the ensuing calculations, which are of course similar to those in the proof of
 Proposition \ref{invt0}.

  Let $f^0 U_{\g_0} , f^1 U_{\g_1} , f^2 U_{\g_2}  \in \Ac$ be such that
 $$
 w(f^0) + w(f^1) + w(f^2) = m \qquad \text{and} \qquad
 \g_0  \g_1  \g_2 =1 \, .
 $$
 Then
 \begin{eqnarray} \nonumber
&\ &b \t (f^0 U_{\g_0} , f^1 U_{\g_1} , f^2 U_{\g_2} )  \, = \cr \cr \nonumber
&=&\t (f^0 f^1| {\g_0}^{-1}  U_{\g_0 \g_1} , f^2 U_{\g_2} ) -
  \t ( f^0 U_{\g_0} , f^1 f^2 | {\g_1}^{-1}  U_{\g_1 \g_2} )  \cr \cr  \nonumber
&+& {\g_2}^{-1} \t (f^2 f^0 | {\g_2}^{-1}  U_{\g_2 \g_0} , f^1 U_{\g_1} )  \cr \cr \nonumber
&=&\int_{z_0}^{\g_0 \g_1 \, z_0} F f^0 f^1| {\g_0}^{-1}  f^2 | {\g_1}^{-1} {\g_0}^{-1} dz
-  \int_{z_0}^{\g_0  \, z_0} F f^0 f^1| {\g_0}^{-1}  f^2 | {\g_1}^{-1} {\g_0}^{-1} dz
 \cr \cr \nonumber
 &+&{\g_2}^{-1}
  \int_{z_0}^{\g_2 \g_0 \, z_0} F f^2 f^0 | {\g_2}^{-1}  f^1 | {\g_0}^{-1} {\g_2}^{-1} dz
 \cr \cr \nonumber
 &=&\int_{\g_0 z_0}^{\g_0 \g_1 z_0} F f^0 f^1| {\g_0}^{-1}  f^2 | {\g_2} dz
 +  \int_{z_0}^{\g_2 \g_0 \, z_0} F | {\g_2}^{-1}
  f^2 f^0 | {\g_2}^{-1}  f^1 | {\g_0}^{-1} {\g_2}^{-1} dz  \cr \cr \nonumber
  &=&\int_{\g_0 \, z_0}^{{\g_2}^{-1} \, z_0} F f^0 f^1| {\g_0}^{-1}  f^2 | {\g_2} dz
  +  \int_{{\g_2}^{-1} \, z_0}^{ \g_0 \, z_0} F
    f^2 | \g_2  f^0  f^1 | {\g_0}^{-1} dz  \quad = \quad 0 \, ,
   \end{eqnarray}
 and so $\, b \t = 0 $.
It is also cyclic, because for $f^0 U_{\g_0} , f^1 U_{\g_1}  \in \Ac$ such that
\begin{equation}  \label{m}
 w(f^0) + w(f^1) = m \qquad \text{and} \qquad
 \g_0  \g_1 =1
 \end{equation}
 one has
  \begin{eqnarray} \nonumber
 &\ & \lb_1 \t (f^0 U_{\g_0} , f^1 U_{\g_1} ) =
  -  {\g_1}^{-1} \t (f^1 U_{\g_1} , f^0 U_{\g_0} ) = \cr \cr \nonumber
  &=&  -  {\g_1}^{-1} \int_{z_0}^{ \g_1 \, z_0} F f^1 f^0 | {\g_1}^{-1}  dz
 = -  \int_{z_0}^{ \g_1 \, z_0} F | {\g_1}^{-1} f^1 f^0 | {\g_1}^{-1}  dz  \cr \cr \nonumber
&=&    \int_{z_0}^{ \g_0 \, z_0}   F  f^1 | {\g_1} f^0   dz
\,   = \, \t (f^0 U_{\g_0} , f^1 U_{\g_1} ) \, .
  \end{eqnarray}
 \end{proof}

\bigskip

\subsection{Cyclic $1$-cocycles with coefficients: base point at cusps}

\medskip

We now extend the construction
of subsection \ref{cusp}, allowing the
base point $z_0$ to belong to the `arithmetic' boundary
of the upper half plane $P^1(\Qb)$,
to the general case of weight $m \geq 2$. To this end
we shall just apply the same procedure to
the cocycle $\t = \t_{m-2}$.
For  $f^0 U_{\g_0} , f^1 U_{\g_1}  \in \Ac$ such that
 \begin{equation}  \label{m'}
  w(f^0) + w(f^1) = m \qquad \text{and} \qquad
 \g_0  \g_1 =1 \
 \end{equation}
it is given by the integral
\begin{equation} \nonumber
\t (f^0 U_{\g_0} , f^1 U_{\g_1} ) \, = \, \int_{z_0}^{ \g_0 \, z_0}  F
 f^0 \, f^1 | {\g_1}  \, dz \, ,
 \end{equation}
 where $F=F_{m-2}$.
 Taking its derivative with respect to $z_0 \in \Hb$ gives
 \begin{eqnarray}   \nonumber
\frac{d}{d z_0} \t (f^0 U_{\g_0} , f^1 U_{\g_1}) &=&  F | {\g_0}  (z_0) \,
 (f^0 | {\g_0} \,  f^1)(z_0) - F(z_0) \, (f^0 \, f^1 | {\g_1} ) (z_0) \cr \cr
 &=&  {\g_0} \cdot F(z_0) \,  (f^0 | {\g_0} \,  f^1)(z_0) -
 F(z_0) \, (f^0 \, f^1 | {\g_1} ) (z_0) \cr \cr
 &=&   - b \e_F \, (f^0 U_{\g_0} , f^1 U_{\g_1} ) ,
 \end{eqnarray}
where
 \begin{equation}   \nonumber
\e_F \,  = \, F(z_0) \, \e \, .
 \end{equation}
As above, we split it into two functionals
 \begin{equation}   \nonumber
\e_F \,  = \,  F(z_0) \, {\ab}_0 \, + \,  F(z_0) \, \wt{\e}  \, ,
 \end{equation}
and regularize the cocycle $\t$ by adding the coboundaries of suitable
anti-derivatives of the two components:
\begin{equation} \nonumber
\wt{\t} \,  = \,  \t  \, +    \, b\left(\check{F} (z_0) \, {\ab}_0 \right)\, - \,
 \int_{z_0}^{i \ify} b \left( F(z) \,   \wt{\e} \right) \, dz \, ,
 \end{equation}
where
 \begin{equation}   \label{tF}
\check{F}_{m-2} (z_0)  \,  = \,  \int_0^{z_0}  F(z) dz \, = \,
\sum_{k=0}^{m-2} \frac{(m-2)!}{(k+1)! (m-k-2)!} \, z_0^{k+1} \, T_1^k \, T_2^{m-k-2} \, .
 \end{equation}
Explicitly, for  $f^0 U_{\g_0} , f^1 U_{\g_1}  \in \Ac$ satisfying (\ref{m'}),
\begin{eqnarray} \label{ttm}
&\ & \wt{\t}_{m-2} (f^0 U_{\g_0} , f^1 U_{\g_1} ) = \cr \cr
&=&  \int_{z_0}^{ \g_0 \, z_0} F
 f^0\,  f^1 | {\g_1} \, dz \, + \,
\check{F}(z_0) \, {\ab}_0 ( f^0 \, f^1 | {\g_1})
  - \g_0\cdot \check{F}(z_0)  \, {\ab}_0 (f^0 | {\g_0} \,  f^1) \cr \cr
 &-& \int_{z_0}^{i \ify} F \, \wt{\e} (f^0 \, f^1 | {\g_1}) \, dz \,
  +  \, \int_{z_0}^{i \ify}  \g_0\cdot F \,  \wt{\e} (f^0 | {\g_0} \,  f^1) \, dz  \, .
 \end{eqnarray}
 Since by its very definition $ \wt{\t}_{m-2}  \in C^1 (\Ac , W_{m-2}) \,$ differs
 from the cocycle $ \t_{m-2}  \in ZC^1 (\Ac , W_{m-2}) \,$ by a coboundary, it is
 itself a cyclic cocycle.
 \bigskip

\section{Transgression for degenerate actions}

The higher weight counterparts of the above results involve
cyclic cocycles with
coefficients that are invariant under `degenerate' actions of $\Hc_1$, associated
to modular forms of arbitrary weight. In the case of cusp forms,
the corresponding
cup products by the universal Godbillon-Vey
cocycle $\d_1 \in ZC^1 (\Hc_1 ; \d, 1)$ transgress to
 $1$-dimensional cohomology classes of congruence subgroups,
 implementing the Eichler-Shimura isomorphism, while
the degenerate actions corresponding to Eisenstein series give rise
to generalized functions of Rademacher type (cf. \cite{Nak}) as well
as to the Eisenstein cocycle of Stevens \cite{Ste}.

\subsection{Degenerate actions of $\Hc_1$ of higher weight}

We recall a few basic facts about cocycle perturbations of
Hopf actions. Given a Hopf algebra $\Hc$ and an algebra
 $\Ac$ endowed with a Hopf action of $\Hc$,  a  $1$-cocycle
 $ \, u \in Z^1 ( \Hc , \Ac ) \, $ is an
invertible element of the convolution algebra $ \,  \Lc( \Hc ,\Ac) \, $ of linear maps from
$\Hc$ to $\Ac$, satisfying
\begin{equation} \label{1coc}
u(h \, h') = \sum u(h_{(1)}) \, h_{(2)}(u(h')) \, , \quad \forall h
\in \Hc \, .
\end{equation}
\noindent The conjugate under $ \, u \in Z^1 ( \Hc , \Ac ) \, $ of the
original action of $\Hc$ on $\Ac$ is given by
\begin{equation} \label{conj}
\tilde{h}(a):= \sum u(h_{(1)}) \, h_{(2)} (a) \, u^{-1}(h_{(3)})
\end{equation}
\medskip

The standard action of $\Hc_1 $ commutes
with the natural coaction of $G^+ (\mathbb{Q})$ on $\Ac_{G^+
(\mathbb{Q})}$. It is natural to only consider cocycles $u$  with the
same property. The values of
such a $1$-cocycle $u$ on generators must belong to the subalgebra
$\Mc \subset \Ac_{G^+ (\mathbb{Q})}$, and  have to be of the
following form:
$$
u(X)= \th \in \Mc_2, \;u(Y)= \lambda \in \mathbb{C}, \; u(\delta_1)=
\om \in \Mc_2 \, .
$$
By \cite[Prop. 11]{CMmha}, for each such data
 $\th  \in \Mc_2, \; \lambda \in \mathbb{C}, \;
 \om \in \Mc_2 \,$,
\begin{itemize}
\item[$1^{0}$] there exists a unique $1$-cocycle
$\displaystyle u \in \Lc( \Hc_1 ,\Ac_{G^+ (\mathbb{Q})})$ such
that
$$
u(X)= \th , \;u(Y)= \lambda , \; u(\delta_1)= \om  \, ;
$$
\item[$2^{0}$] the conjugate under $u$ of the
       action of $\Hc_1$ is given
       on generators as follows:
\begin{eqnarray}
\tilde{Y} = Y \, , \quad \tilde{X}(a) &=& X(a)+[(\th - \lambda \om),a]
- \lambda\, \delta_1(a)+ \om Y(a) \, ,  \nonumber \\
 \tilde{\delta_1}(a) &=& \,\delta_1(a) \label{diff}
+ [\om,a] \, , \qquad a \in \Ac_{G^+ (\mathbb{Q})}\, ;
\end{eqnarray}
\item[$3^{0}$] The conjugate under $u$  of $\delta'_2$ is given
by the operator
$$
\tilde{\delta'_2}(a) = [X(\om) +\frac{ \om^2}{2}- \, \Om_4,\; a] \, ,
\qquad a \in \Ac_{G^+ (\mathbb{Q})} $$ and there is no choice of
$u$ for which $\tilde{\delta'_2} = 0$.
\end{itemize}
\medskip

\noindent The actions described above were called
  \textit{projective} because $\, \tilde{\d}'_2 $ acts
by an inner transformation.
\bigskip

For our purposes, it is the difference ${\wt \d_1} - \d_1$ which matters.
In view of (\ref{diff}), we may as well start from the trivial action
and may also assume $\th  = 0$, $\lb = 0$.
\medskip

We call the \textit{trivial action of
weight} $m$ of $\Hc_1$ on $\, \Ac_{G^+ (\mathbb{Q})}  \,$ the
action defined by
\begin{equation}  \label{aff0}
Y (a) \, = \, \frac{w(f)}{m} \, a ,
\qquad X (a) \, = \, 0 \, , \qquad
\delta_n (a) \, = \, 0 , \quad n \geq 1 \,
\end{equation}
\medskip

To any modular form $\, \om \in {\Mc}_m \,$ we shall now associate
a `degenerate' action of the Hopf algebra $\Hc_1$ on the crossed product algebra
$\, \Ac_{G^+ (\mathbb{Q})}  \,$ as follows.
\bigskip

\begin{proposition} \label{affac}
Let $\, \om \in {\Mc} \,$ be a modular of weight $\, w(\om) \, = \, m $.

\begin{itemize}
\item[$1^{0}$.]  There exists a unique $1$-cocycle
$\displaystyle \, u = u_{\om} \in Z^1 ( \Hc_1 ,\Ac_{G^+ (\mathbb{Q})}) \, $ such that
\begin{equation*}
u(Y) \, = \, 0 , \qquad u(X) \, = \, 0 , \qquad u(\delta_1) \, = \, \om  \, .
\end{equation*}

\item[$2^{0}$.] The conjugate under $\, u_{\om} \,$ of the trivial action of weight $m$
of $\Hc_1$  on $\, \Ac = \Ac_{G^+ (\mathbb{Q})}  \,$ is determined by
\begin{equation} \label{aff1}
\tilde{Y} (a) = \frac{w(a)}{m} \, a ,
\quad \tilde{X} (a) = \om \, \tilde{Y} (a), \quad
 \tilde{\delta_1}(a) =  [\om , a], \quad \fl \, a \in \Ac \, .
\end{equation}
Under this action, one has, for any $\, n \geq 1 \,$,
\begin{eqnarray}
\tilde{\delta_n}(a) &=& (n-1)! \, \om^{n-1} \, [\om , a] , \qquad \fl \, a \in \Ac \, , \label{aff2} \\
\text{resp.} \quad
\tilde{\delta_n}(f \, U_{\g}^*) &=&  X^{n-1} (\om \, - \, \om \vert \g) \, f \, U_{\g}^*  ,
\quad \fl \, f \, U_{\g}^* \in \Ac \, . \label{aff3}
\end{eqnarray}
\end{itemize}
\end{proposition}

\medskip
\begin{proof}
The first statement is a variant of \cite[Proposition 11, $1^0$]{CMmha}.
Its proof gives the following
explicit description of the $1$-cocycle
$\displaystyle \, u \in Z^1 ( \Hc_1 ,\Ac_{G^+ (\mathbb{Q})}) $.

\noindent First, recall that any element of $\Hc_1$ can be uniquely written as a
linear combination of monomials of the form $\displaystyle
P(\d_1,\d_2,...,\d_{\ell}) \, X^n \, Y^m \,$. We then let
$\om^{(k)}$ be defined by
induction by
\begin{equation*}
\om^{(1)} \, := \, \om , \qquad \om^{(k+1)} \, :=\, \om\, Y(\om^{(k)})
\, =\, k! \, \om^{k+1}  ,
\quad \fl \, k \geq 1 \, .
\end{equation*}
Then
\begin{equation*}
u \left( P(\d_1,\d_2,...,\d_{\ell})\, X^n \, Y^m \right) \, =\,
\begin{cases}
P(\om^{(1)},\om^{(2)},...,\om^{({\ell})}) \quad &\text{ if }  \quad n=m=0 , \\
 0  \quad &\text{ otherwise} \,  ,
\end{cases}
\end{equation*}
while its inverse $\, u^{-1} \,$  is
given by
\begin{equation*}
u^{-1}\left( P(\d_1,\d_2,...,\d_{\ell})\, X^n \, Y^m \right) \,=\,
\begin{cases}
P( - \om, 0,...,0) \quad &\text{ if }  \quad n=m=0 , \\
 0  \quad &\text{ otherwise} \,  .
\end{cases}
\end{equation*}
\smallskip

The formulae (\ref{aff1}) are obtained from $1^0$  above and the definitions
(\ref{conj}), (\ref{aff0})
(cf. also \cite[Proposition 11, $2^0$]{CMmha}).
 Finally (\ref{aff2}),
and hence (\ref{aff3}), follows by straightforward computation.
\end{proof}
\medskip

\noindent The actions as in $2^0$  above are \textit{degenerate}, in the sense
 that $\, \tilde{\d}_1 \,$ acts by an inner transformation.
 \bigskip

 \begin{lemma} \label{invtm}
 With respect to the degenerate action of weight $m$, one has
  \begin{equation} \label{invp}
 \t (h_{(1)}(a^0) , h_{(2)}(a^1)) \, = \, \d(h) \, \t (a^0 , a^1) , \qquad \fl \, h \in \Hc_1 , \,  \,
 a^0 , a^1 \in \Ac \, ,
 \end{equation}
 for $\t = \t_{m-2}$ or $\t = \wt{\t}_{m-2}$.
 \end{lemma}
 \smallskip

 \begin{proof} We shall check the $\Hc_1$-invariance
property (\ref{invp}) for  $\t = \t_{m-2}$.
 As in the
proof of Proposition \ref{invt0},  it suffices to verify it on the generators.

Starting with ${\wt Y}$,  for
$f^0 U_{\g_0} , f^1 U_{\g_1} $ satisfying (\ref{m})
one has
 \begin{eqnarray} \nonumber
&\ &  \t ({\wt Y}(f^0 U_{\g_0}) , f^1 U_{\g_1} ) +
 \t (f^0 U_{\g_0} , {\wt Y}(f^1 U_{\g_1} )) =
\cr \cr \nonumber
&=&  \int_{z_0}^{ \g_0 \, z_0} F  {\wt Y}( f^0 f^1 | {\g_0}^{-1})  dz
= \frac{w(f^0) + w(f^1)}{m} \int_{z_0}^{ \g_0 \, z_0} F  f^0 f^1 | {\g_0}^{-1}  dz
\cr \cr \nonumber
&=& \int_{z_0}^{ \g_0 \, z_0} F  f^0 f^1 | {\g_0}^{-1}  dz \quad = \quad
\d(Y) \, \t (f^0 U_{\g_0} , f^1U_{\g_1} ) \, .
\end{eqnarray}
 Passing to ${\wt X}$, the identity  (\ref{invp}) is nontrivial only if
   $f^0 U_{\g_0} , f^1 U_{\g_1}  \in \Ac$ satisfy
 \begin{equation}  \label{0}
 w(f^0) + w(f^1) = 0 \qquad \text{and} \qquad
 \g_0  \g_1 =1 \, .
 \end{equation}
 One then has
 \begin{eqnarray} \nonumber
&\ & \t ({\wt X}(f^0) U_{\g_0} , f^1 U_{\g_1} ) +  \t (f^0 U_{\g_0} , {\wt X}(f^1) U_{\g_1} )
+  \t (\d_1(f^0 U_{\g_0}) , {\wt Y}(f^1) U_{\g_1} ) =
\cr \cr \nonumber
&=&  \int_{z_0}^{ \g_0 \, z_0} F \, \om {\wt Y}( f^0) f^1 | {\g_0}^{-1}  dz +
 \int_{z_0}^{ \g_0 \, z_0} F  \, f^0 \om | {\g_0}^{-1} {\wt Y}( f^1) | {\g_0}^{-1}  dz
 \cr \cr \nonumber
&+&  \int_{z_0}^{ \g_0 \, z_0} F  \, (\om - \om | {\g_0}^{-1})
f^0 {\wt Y}( f^1) | {\g_0}^{-1}  dz \, = \,
\int_{z_0}^{ \g_0 \, z_0} F  \, \om
 {\wt Y}( f^0 \, f^1 | {\g_0}^{-1})  dz   \cr \cr \nonumber
&=&  \frac{w(f^0) + w(f^1)}{m} \int_{z_0}^{ \g_0 \, z_0} F  f^0 f^1 | {\g_0}^{-1}  dz
\quad = \quad 0 \, ,
 \end{eqnarray}
the vanishing being a consequence of (\ref{0}).

 Finally, in the case of ${\wt \d_1}$ and
 with  $f^0 U_{\g_0} , f^1 U_{\g_1} $ as in (\ref{0}),
  if the action is degenerate, one has
  \begin{eqnarray} \nonumber
&\ &  \t ({\wt \d_1}(f^0 U_{\g_0}) , f^1 U_{\g_1} ) +
 \t (f^0 U_{\g_0} , {\wt \d_1}(f^1 U_{\g_1} )) =
\cr \cr \nonumber
&=&  \int_{z_0}^{ \g_0 \, z_0} F  (\om - \om |  {\g_1} )   f^0 f^1 | {\g_1}  dz
 +  \int_{z_0}^{ \g_0 \, z_0} F  f^0 (\om - \om |  {\g_0} ) |  {\g_1}   f^1 | {\g_1}  dz
\, = \, 0 \, .
 \end{eqnarray}
 \end{proof}
\bigskip

\begin{remark} \label{tA}
{\rm The fact that the only modular forms of weight $0$ are
the constants was used just in the case of the standard
action. Thus, for all the degenerate actions, the statement remains valid for
the enlarged crossed product $\, {\wt \Ac} = {\wt \Mc} \ltimes
\GL^+(2, \Qb)$, where ${\wt \Mc} $ is the algebra consisting of  meromorphic
 modular forms of all levels whose poles are located on cusps.}
 \end{remark}

\subsection{Transgression in the higher weight case}

We now consider the degenerate action associated to
a modular form $\om \in \Mc_m$ with $m \geq 2$
and define the cup product of
$\d_1$ with  the invariant $1$-trace
$\t $ as in Lemma \ref{invtm} by the formula
\begin{equation} \nonumber
gv_{\om} \, (a^0 , a^1 , a^2)  =  \t (a^0 \, \wt{\d_1} (a^1), \,a^2)
=  \t (a^0 \, \om a^1- a^1 \om, \,a^2) \, ,
\quad a^0 , a^1 , a^2 \in \Ac \,
 \end{equation}
\medskip

\begin{lemma}  \label{lgvm}
One has $gv_{\om} \in ZC^2 (\Ac, W_{m-2})$, that is
$$ \nonumber
\quad b (gv_{\om}  )\, = \, 0 \quad \text{and} \quad B (gv_{\om} ) \, = \, 0 \, .
$$
\end{lemma}
\smallskip

\begin{proof} As in the proof of Lemma \ref{lgv},
using the primitivity of $\d_1$, one has
\begin{eqnarray} \nonumber
&\ &b (gv_{\om}  )\,(a^0 , a^1 , a^2 , a^3) \, =\, \t \, (a^0 a^1 \d_1 (a^2) , a^3 ) \, - \,
\t \, (a^0  \d_1 (a^1 a^2) , a^3 ) \cr \cr
&+& \t \, (a^0 \d_1 (a^1),  a^2 a^3 ) \, - \,
S(a_{(-1)}^3)\,  \t \, (a_{(0)}^3 a^0 \d_1(a^1) , a^2 )  \cr \cr
&=& - \t \, (a^0 \d_1 (a^1)  a^2 , a^3 ) \, + \,
  \t \, (a^0  \d_1 (a^1),  a^2 a^3 ) \cr \cr
 &-&S(a_{(-1)}^3)\, \t \, (a_{(0)}^3 a^0 \d_1(a^1) , a^2 )
\, = \,- b\t \, (a^0 \d_1 (a^1),  a^2 , a^3 ) \, = \, 0 \, .
\end{eqnarray}

Next, since $\, \t (a, 1) \, = \, 0$ for any $a\in \Ac$, one has
\begin{eqnarray} \nonumber
B (gv_{\om}  )\,(a^0 , a^1) &=& g_{\om} v \, (1, a^0, a^1) \, - \,
S(a_{(-1)}^1)\,gv_{\om}  \, (1, a_{(0)}^1, a^0) \cr \cr
 &=& \t \, ( \d_1 (a^0) , a^1 ) \, - \,
S(a_{(-1)}^1)\, \t \, (\d_1 (a_{(0)}^1) , a^0 ) \, .
 \end{eqnarray}
 We now use the fact that the action of $\d_1$ commutes with the
 coaction of $\GL^+ (2, \Qb)$, to continue
 \begin{eqnarray} \nonumber
B (gv_{\om}  )\,(a^0 , a^1) &=&  \t \, ( \d_1 (a^0) , a^1 ) \, - \,
S(\d_1(a^1)_{(-1)}) \, \t \, (\d_1 (a^1)_{(0)} , a^0 ) \cr \cr
 &=& \t \, ( \d_1 (a^0) , a^1 ) \, + \,  \t \, ( a^0 , \d_1 (a^1) ) \, = \, 0 \, ,
 \end{eqnarray}
the vanishing being a consequence of the $\Hc_1$-invariance
property of  $\t \in ZC^1 (\Ac, W_{m-2})$.
\end{proof}
\bigskip

Because degenerate actions are perturbations
of the trivial action of weight $m$, one expects the cup product
$$ \nonumber
[gv_{\om} ] \, = \, [\d_1] \# [\t_{m-2} ]
$$
to vanish. This is indeed the case, the
vanishing being a consequence of a general transgression formula
for degenerate actions. To state it, we introduce
the cochain $tgv_{\om} \in C^1 (\Ac, W_{m-2})$,
\begin{equation} \label{tgv}
 tgv_{\om} (a^0 , a^1) \, = \, - \t (a^0 ,  \om \, a^1) \, , \qquad a^0 , a^1  \in \Ac \,
\end{equation}
\medskip

\begin{proposition} \label{traff}
For any $\om \in \Mc_m$,
 $$
  gv_{\om} \,   = \, b\, (tgv_{\om}) \quad \text{and} \quad B\, (tgv_{\om})\, =\, 0 \, .
 $$
  \end{proposition}

\begin{proof} One has
\begin{equation} \nonumber
 b (tgv_{\om}) (a^0 , a^1 , a^2) = - \t (a^0  a^1 , \, \om  a^2) + \t (a^0  , \, \om a^1 a^2)
 - a_{(-1)}^2  \t (a^2  a^0 , \, \om  a^1) ;
 \end{equation}
 after replacing the first term in the right hand side using
 \begin{eqnarray} \nonumber
 0\, &=& \, b\t (a^0 ,  a^1 ,  \om  a^2) = \cr \cr
 &=& \t (a^0  a^1 , \, \om  a^2) - \t (a^0 , \, a^1 \om  a^2)
 + a_{(-1)}^2  \t (\om a^2  a^0 , \, a^1) ,
 \end{eqnarray}
one obtains
\begin{eqnarray} \nonumber
&\ & b (tgv_{\om}) \, (a^0 , a^1 , a^2) \, = \,   \t (a^0  , \, \om a^1 a^2) - \t (a^0 , \, a^1 \om  a^2)
\,  + \cr \cr
 &+& a_{(-1)}^2 \left( \t (\om a^2  a^0 , \, a^1) -
 \t (a^2  a^0 , \, \om  a^1) \right) \cr \cr
 &=& \t (a^0  , \, \wt{\d_1}( a^1) a^2) + a_{(-1)}^2 \left( \t (\om a^2  a^0 , \, a^1) -
 \t (a^2  a^0 , \, \om  a^1) \right) .
 \end{eqnarray}
We now use the identity
 \begin{equation} \nonumber
 0= b\t  (a^0  , \wt{\d_1}( a^1),  a^2)  = \t  (a^0 \wt{\d_1}( a^1),  a^2) -
 \t  (a^0  , \wt{\d_1}( a^1) a^2) + a_{(-1)}^2 \t  ( a^2 a^0  , \wt{\d_1}( a^1))
  \end{equation}
to replace the term $\t (a^0  , \, \wt{\d_1}( a^1) a^2)$, thus obtaining
\begin{eqnarray} \label{btgv}
&\ &b (tgv_{\om}) \, (a^0 , a^1 , a^2)=  \t  (a^0 \wt{\d_1}( a^1),  a^2) \, +
 \cr \cr
 &+& a_{(-1)}^2 \left(\t  ( a^2 a^0  , \wt{\d_1}( a^1)) +  \t (\om a^2  a^0 , \, a^1) -
 \t (a^2  a^0 , \, \om  a^1) \right) \cr \cr
&=&  \t  (a^0 \wt{\d_1}( a^1),  a^2)  +
a_{(-1)}^2 \left(-\t  ( a^2 a^0  , a^1\om) +  \t (\om a^2  a^0 , \, a^1) \right) .
   \end{eqnarray}
Next, using
  \begin{eqnarray} \nonumber
 0= b\t  (c^0  ,  c^1,  \om) &=& \t  (c^0 c^1,  \om) -
 \t  (c^0  , c^1\om) +  \t  ( \om c^0  , c^1) \cr \cr
 &=& - \t  (c^0  , c^1\om) +  \t  ( \om c^0  , c^1) \, ,
  \end{eqnarray}
  for $c^0 = a^2 a^0$ and $c^1 = a^1$, one sees that the term in
  paranthesis vanishes, and therefore (\ref{btgv}) reduces to
  \begin{equation} \nonumber
 b (tgv_{\om}) \, (a^0 , a^1 , a^2)\, = \, \t  (a^0 \wt{\d_1}( a^1),  a^2)
 \, = \, gv_{\om} \,  (a^0 , a^1 , a^2) \, .
  \end{equation}

 On the other hand, by the very definition of $\t =
 \t_{m-2}$, one has
   \begin{equation} \nonumber
 B (tgv_{\om}) \, (a^0 )\, = \, -  \t (1, \om a^0) \, - \t (a^0 , \om) \, = \, 0 \, .
  \end{equation}
\end{proof}
\bigskip

Let $\om \in \Mc_m (\G)$, for some congruence subgroup $\G \sbs \SL(2, \Zb)$.
The transgression proper occurs when  $gv_{\om}$ is restricted to the
subalgebra
$$
\Ac^{\G} \, = \, \Mc \ltimes \G  \, ,
$$
on which $\wt{\d_1} $ vanishes. Then $ tgv_{\om}$ becomes a cocycle  in
$ZC^2 (\Ac^{\G}, W_{m-2})$ and its restriction to
$ \Ac^{\G}_0 = \Cb [\G] $
\begin{equation} \nonumber
 tgv_{\om}  (U_{\g_0} , U_{\g_1}) =
 - \t (U_{\g_0} , \om \, U_{\g_1}) =
- \int_{z_0}^{\g_0 z_0} F \om \, dz \, , \qquad \g_0 \g_1 = 1,
\end{equation}
 gives a group cocycle in  $Z^1 (\G, W_{m-2})$. After a change of sign,
 we denote this cocycle
 \begin{equation} \label{TGV}
 TGV_{\om} (\g) := \, \int_{z_0}^{\g z_0} F \om \, dz \, , \qquad \fl \,
 \g \in \G .
 \end{equation}
 Taking its real part, one obtains the linear map
 \begin{equation} \label{ltgv}
 \om \in  \Mc_m (\G) \longmapsto
 ES(\om) := [\Re TGV_{\om}] \in H^1 (\G, W_{m-2} (\Rb)) \, ,
 \end{equation}
 whose restriction to the cuspidal subspace $\Mc^0_m (\G) $ gives
 the Eichler-Shimura embedding of
  $  \Mc^0_m (\G)$ into  $H^1 (\G, W_{m-2} (\Rb))$.
   \bigskip

 In the remainder of this section we shall look closer at the restriction
 of the assignment (\ref{ltgv}) to
 the Eisenstein subspace $ \Ec_m (\G) \sbs  \Mc_m (\G)$. To this end,
 we proceed to recall the construction of the Hecke lattice for higher weight
  Eisenstein series.

For $\, {\ab} = (a_1, a_2) \in ( \mathbb Q \slash
\mathbb Z )^2  $, the holomorphic
Eisenstein series $G^{(m)}_{\ab}$ of  weight $m > 2$
is defined by the absolutely convergent series
\begin{equation*}
 G^{(m)}_{\ab} (z) := \, \sum_{ {\kb} \in {\mathbb Q}^2 \setminus {\0b},
\, {\kb} \equiv {\ab} \, (\text{mod} \,1) }
 (k_{1} z + k_{2})^{-m}  \, ;
\end{equation*}
averaging over $ \left(\frac{1}{N} \Zb \slash \Zb\right)^2$, and using as
 weights the additive characters  (see (\ref{chi}))
 $\, \left\{\chi_{\xb} \, ;
 \, {\xb}   \in \left( \frac{1}{N} \Zb \slash \Zb \right)^2 \right \} $,
gives rise to the series
\begin{equation} \label{mbasis}
  \phi^{(m)}_{\xb} (z) :=  \frac{(m-1)!}{(2 \pi i N)^m} \sum_{{\ab} \in
 \left( \frac{1}{N} \mathbb Z \slash \mathbb Z \right)^2} \,
 \chi_{\xb} ({\ab}) \cdot  G^{(m)}_{\ab} (z) \, .
\end{equation}
\medskip

  We now apply Proposition \ref{traff} for
 the case of the degenerate action
 defined by an Eisenstein series  $\, \phi^{(m)}_{\xb} $,
  $\, {\xb} \in \left( \frac{1}{N} \Zb \slash \Zb \right)^2 $,
 and with respect to the $\Hc_1$-invariant cocycle
 $\wt{\t} = \wt{\t}_{m-2}$ (see (\ref{ttm})). Restricting
 to $\Cb [\G (N)]$ one obtains the transgressed group cocycle
 $ \wt{TGV}_{\phi_{\xb}} \in Z^1 (\G (N), W_{m-2})$.
 We then take its real part
 $\Phi^{(m)}_{\xb} := \Re
 \wt{TGV}_{\phi_{\xb}} \in Z^1 (\G (N), W_{m-2} (\Rb))$, which still satisfies the cocycle relation
  \begin{equation} \label{gc}
\Phi^{(m)}_{\xb} (\a \b) \, = \,  \Phi^{(m)}_{\xb} (\a ) \, + \, \a  \cdot \Phi^{(m)}_{\xb} (\b) ,
\qquad  \a , \b \in \G(N) \, .
 \end{equation}
 We shall show that  $\Phi^{(m)}_{\xb} $
  coincides with the
 generalized Rademacher function of \cite[\S 2]{Nak}, and
 in particular $\Phi^{(m)}_{\xb}  \in Z^1 (\G (N), W_{m-2} (\Qb))$.
  As a matter of fact, we shall explicitly compute
  its extension to $ \GL^+(2, \Qb) \,$,
$$
\Phi^{(m)}_{\xb} \, = \, \Re \Psi^{(m)}_{\xb} \, .
$$
In turn,
$$\Psi^{(m)}_{\xb} (\g) \, = \, \wt{TGV}_{\phi_{\xb}} \, ,
\qquad \g \in \GL^+(2, \Qb) \, ,
$$
is
obtained by specializing the formula (\ref{ttm}) to $f^0 = \phi^{(m)}_{\xb}$, $f^1 = 1$
and thus is given by the expression
 \begin{eqnarray} \label{Pix}
 \Psi^{(m)}_{\xb} (\g) &=& \int_{z_0}^{\g z_0} F (z) \, \phi^{(m)}_{\xb}(z) dz   +
 \check{F}(z_0) \, {\ab}_0 (\phi^{(m)}_{\xb})
  -   \g \cdot \check{F}(z_0) \, {\ab}_0 (\phi^{(m)}_{\xb} | \g) \cr \cr
 &-& \int_{z_0}^{i \ify}  F(z) \, \wt{\e}({\phi}^{(m)}_{\xb})(z) dz + \int_{z_0}^{i \ify}
   \g \cdot F (z)  \,  \wt{\e}({\phi}^{(m)}_{\xb} | \g)(z) dz  \, ,
 \end{eqnarray}
 which is independent of $z_0 \in \Hb$.
 The coboundary relation in Proposition \ref{traff} is
 equivalent with the following
 cocycle property
   \begin{equation} \label{tgc}
\Psi^{(m)}_{\xb} (\a \b) \, = \,  \Psi^{(m)}_{\xb} (\a ) \, + \, \a  \cdot \Psi^{(m)}_{\xb | \a} (\b) ,
\qquad  \a , \b \in \GL^+(2, \Qb) \, ,
 \end{equation}
 and the collection
 $\{ \Phi^{(m)}_{\xb} \, ; \, {\xb} \in \Qb^2/\Zb^2 \}$
 is equivalent to the distribution valued
 \textit{Eisenstein cocycle} of Stevens \cite{Ste}.
\bigskip

 To state the precise result, we need to recall one
 more definition, that of a
\textit{generalized higher Rademacher-Dedekind sum}
(cf.~\cite{RG, HWZ, Ste, Nak});
given $\, a, c \in \Zb\, $ with $(a, c) = 1$ and $c \geq 1$, $ 0 < k < m \in \Nb$ and
$ {\xb}   \in \left( \frac{1}{N} \Zb \slash \Zb \right)^2  $, it is defined by the expression
\begin{equation} \label{GRD}
  S^{(m-k, k)}_{\xb} (\frac{a}{c})  = \sum_{r=0}^{c-1}
   \frac{{\Bb}_{m-k} \left(\frac{x_1+r}{c}\right)}{m-k} \,
    \frac{{\Bb}_{k} \left(x_2 + a \frac{x_1 + r}{c} \right)}{k}  \, ,
   \end{equation}
where $ {\Bb}_j : \Rb \ra \Rb$ is the $j$-th periodic Bernoulli function
 $$
 {\Bb}_j (x) \,  = \, B_j (x - [x]) \, , \qquad x \in \Rb \, .
 $$
\medskip

\begin{theorem} \label{RadFi}
Let $m \geq 2$ and let ${\xb} \in \Qb^2/\Zb^2 $, with $\, {\xb} \neq 0$
if $\, m=2$.

$1^0$. For $\b =  \begin{pmatrix} a &b\\
 0 &d  \end{pmatrix} \in B (\Qb)$,
 \begin{equation} \nonumber
   \Phi^{(m)}_{\xb} (\b) \, = \, - \frac{ {\Bb}_m  (x_1)}{m} \,
   \int_{0}^{ \frac{b}{d} } (t T_1 + T_2)^{m-2} \,dt  \, .
 \end{equation}

$2^0$. For $\s =  \begin{pmatrix} a &b\\
 c &d  \end{pmatrix} \in \SL(2, \Zb)$, with $c > 0$,
  \begin{equation} \nonumber
 \Phi^{(m)}_{\xb} (\s) = \left \{\begin{matrix} &-&
 \displaystyle{\frac{ {\Bb}_m  (x_1)}{m} \,
   \int_{0}^{ \frac{a}{c} } } (t T_1 + T_2)^{m-2} \,dt
   \qquad \qquad \qquad \qquad \qquad \qquad \qquad \cr \cr \
&-&\displaystyle{\frac{ {\Bb}_m  (a x_1 + cx_2)}{m} \,
   \int_{-\frac{d}{c}}^{0 }} \left(t (aT_1 + cT_2) + bT_1 + d T_2 \right)^{m-2} \,dt
   \qquad \qquad  \cr \cr \
&+&\displaystyle{\sum_{k=0}^{m-2}} (-1)^k
 \begin{pmatrix} m-2\\
 k  \end{pmatrix}
  S^{(m-k-1, k+1)}_{\xb} (\frac{a}{c}) T_1^k (a T_1 + cT_2)^{m-k-2} . \qquad \qquad
\end{matrix} \right.
\end{equation}
\end{theorem}
\medskip

\begin{proof}
$1^0$. Restricting $\Psi^{(m)}_{\xb}$ to
$B (\Qb)$, one has
  \begin{eqnarray} \nonumber
 \Psi^{(m)}_{\xb} (\b) &=& \int_{ z_0}^{\b z_0}  F (z)\, \wt{\e}({\phi}^{(m)}_{\xb})(z)  dz
 +  {\ab}_0 ( \phi^{(m)}_{\xb} )  \int_{ z_0}^{\b z_0}  F (z)  dz  \cr \cr
 &+&   \check{F}(z_0) \, {\ab}_0 (\phi^{(m)}_{\xb})
  -   \b \cdot \check{F}(z_0) \, {\ab}_0 (\phi^{(m)}_{\xb} | \b) \cr \cr
 &-& \int_{z_0}^{i \ify}  F(z)\,  \wt{\e}({\phi}^{(m)}_{\xb})(z) \, dz + \int_{z_0}^{i \ify}
   \b \cdot F (z) \,  \wt{\e} ( \phi^{(m)}_{\xb} | \b)(z) dz  \cr \cr
 &=&  \int_{ z_0}^{\b z_0}  F (z)\, \wt{\e}({\phi}^{(m)}_{\xb})(z)  dz
+  {\ab}_0 ( \phi^{(m)}_{\xb} )\, ( \check{F}(\b z_0) - \check{F}(z_0) ) \cr \cr
 &+&   \check{F}(z_0) \, {\ab}_0 (\phi^{(m)}_{\xb})
  -   \b \cdot \check{F}(z_0) \, {\ab}_0 (\phi^{(m)}_{\xb} | \b) \cr \cr
 &-& \int_{z_0}^{i \ify}  F(z)\,  \wt{\e}({\phi}^{(m)}_{\xb})(z) dz + \int_{z_0}^{i \ify}
   \b \cdot F (z) \,  \wt{\e} ( \phi^{(m)}_{\xb} | \b)(z) dz  \cr \cr
   &=&  \int_{ z_0}^{\b z_0}  F (z) \wt{\e}({\phi}^{(m)}_{\xb})(z)  dz
+  {\ab}_0 ( \phi^{(m)}_{\xb} )  \check{F}(\b z_0)
  -  \b \cdot \check{F}(z_0) {\ab}_0 (\phi^{(m)}_{\xb} | \b) \cr \cr
 &-& \int_{z_0}^{i \ify}  F(z)\,  \wt{\e}({\phi}^{(m)}_{\xb}) (z) dz + \int_{z_0}^{i \ify}
   \b \cdot F (z) \,  \wt{\e} ( \phi^{(m)}_{\xb} | \b)(z) \, dz  \, .
 \end{eqnarray}
 When $z_0 \ra i \ify$ all three integrals vanish. For the remaining
 terms we note that, using the weight $m$
 analogue of (\ref{a0b}),
 one has
  \begin{equation} \nonumber
    {\ab}_0 ( \phi^{(m)}_{\xb} ) \check{F}(\b z_0)
  -  \b \cdot \check{F}(z_0) {\ab}_0 (\phi^{(m)}_{\xb} | \b)
 =  {\ab}_0 ( \phi^{(m)}_{\xb} ) \left( \check{F}(\b z_0) - \left(\frac{a}{d}\right)^{\frac{m}{2}}
\b \cdot \check{F}( z_0)\right).
 \end{equation}
 Since
  \begin{equation} \nonumber
    \frac{d}{d z_0}  \left( \check{F}(\b z_0) - \left(\frac{a}{d}\right)^{\frac{m}{2}}
\b \cdot \check{F}( z_0)\right) \, = \, 0 \, ,
 \end{equation}
 as one can easily check employing (\ref{Fcov}),  the expression is
 constant in $z_0$, hence equal to its value at $z_0 = 0$
  \begin{equation} \nonumber
    \check{F}(\b z_0) - \left(\frac{a}{d}\right)^{\frac{m}{2}}
\b \cdot \check{F}( z_0) \, = \, \check{F}( \frac{b}{d}) \, .
 \end{equation}
 From the known Fourier transform of $\phi^{(m)}_{\xb}$
 (cf. \cite{Nak}, formula (2.1)), one reads that
 \begin{equation} \label{a0phix}
   {\ab}_0 ( \phi^{(m)}_{\xb} ) \, = \, - \frac{ {\Bb}_m (x_1)}{m} \, .
 \end{equation}
 It follows that
  \begin{equation} \label{phib}
   \Phi^{(m)}_{\xb} (\b) \, = \, - \frac{ {\Bb}_m  (x_1)}{m} \,
   \int_{0}^{ \frac{b}{d} } (t T_1 + T_2)^{m-2} \,dt  \,
 \end{equation}
 which proves the first statement.
\medskip

$2^0$. Let us now compute the value of
$ \Psi^{(m)}_{\xb}$ at $\s_0 =  \begin{pmatrix} 0 &-1\\
 1 &\, \, \, 0  \end{pmatrix}$ by choosing in (\ref{Pix}) $z_0 = i$, which is a fixed point of $\s_0$.
 One has
 \begin{eqnarray} \nonumber
 \Psi^{(m)}_{\xb} (\s_0) &=&   {\ab}_0 ( \phi^{(m)}_{\xb} )\,  \check{F}(i)
  -   \s_0 \cdot \check{F}(i) \, {\ab}_0 (\phi^{(m)}_{\xb} | \s_0) \cr \cr
 &-& \int_{i}^{i \ify}  F(z)\,  \wt{\e}(\phi^{(m)}_{\xb})(z) \, dz + \int_{i}^{i \ify}
   \s_0 \cdot F (z) \,  \wt{\e} ( \phi^{(m)}_{\xb} | \s_0)(z) \, dz  \, .
 \end{eqnarray}
 This can be related to the Mellin transform of $F \phi^{(m)}_{\xb}$,
   \begin{equation} \nonumber
  D (F \phi^{(m)}_{\xb}, s) \, = \,  \int_{0}^{i \ify}  F(z)\, \wt{\e}(\phi^{(m)}_{\xb})(z)
\, y^{s-1}  \, dz  \, , \quad z = x+iy \, .
 \end{equation}
More precisely, as in the proof of  \cite[Prop. 2.3.3]{St}, one shows that
  \begin{equation} \nonumber
  \Psi^{(m)}_{\xb} (\s_0) \, = \, - D (F \phi^{(m)}_{\xb}, 1)  \, ,
 \end{equation}
 which in turn can be computed, as in \cite[Prop. 2.2.1]{St},
 from the Fourier transform of f $\phi^{(m)}_{\xb}$. Taking the real part
 one obtains, cf. \cite[Thm. 6.9 (a)]{Ste},
 \begin{equation} \label{phir}
  \Phi^{(m)}_{\xb} (\s_0)  = \sum_{k=0}^{m-2} (-1)^k
   \begin{pmatrix} m -2\\
 k  \end{pmatrix}
   \frac{{\Bb}_{m-k-1} (x_1)}{m-k-1} \, \frac{{\Bb}_{k+1} (x_2)}{k+1}  \,
   T_1^k \, T_2^{m-k-2}
 \end{equation}
 \medskip

 In view of the Bruhat decomposition
 $\, \GL^+(2, \Qb)  =  B(\Qb) \cup B(\Qb) \s_0 B(\Qb) $,
  the expressions (\ref{phib}) and (\ref{phir}), together with  the cocycle relation (\ref{tgc}),
  uniquely determine $\Phi^{(m)}_{\xb} $ and allow its explicit calculation (comp.
  \cite[Ch. 2]{St}, \cite[\S 5]{Ste}, where the weight $2$ case is treated in detail).
  This calculation has in fact been done by Nakamura \cite{Nak}, for a
  cocycle defined in a similar manner but using an alternate construction
  of the Eichler-Shimura period integrals.
  Since both  (\ref{phib}) and (\ref{phir})
   agree with Nakamura's formula \cite[(2.8)]{Nak} for the
restriction of $\Phi^{(m)}_{\xb} $ to $\SL(2, \Zb)$,
 one can conclude that $\Phi^{(m)}_{\xb} (\s) $ is given by the cited formula
 for any $\s \in \SL(2, \Zb)$.
 \end{proof}
\bigskip

A similar construction can be performed starting with any
holomorphic
function $F: \Hb \ra W$ that satisfies the covariance law (\ref{Fcov}),
\begin{equation} \nonumber
 F \vert g \, (z) \equiv \det (g)^{- \frac{m-2}{2}}  (cz + d)^{m-2} \, F (g z)
  = g \cdot F  (z) \, , \qquad g \in GL ^+(2, \mathbb{R}).
  \end{equation}
Of particular interest, in view of its connection to special values at non-positive
integers of partial zeta functions for real quadratic fields (see \cite[\S 7]{Ste}),
 is the following choice (cf. \cite[(6.4) (b)]{Ste}).
For  $m = 2n \in 2\Nb$, one takes
 $$
 W_{n-1, n-1} \, = \, (W_{n-1} \ot W_{n-1})^{\rm sym}  \, \sbs \,W_{n-1} \ot W_{n-1} \, ,
 $$
i.e. the subspace fixed by the involution $P_1 \ot P_2 \mapsto P_2 \ot P_1$;
its elements may be identified with the
homogeneous polynomials $P \in \Cb [T_1, T_2, T_3, T_4]$
of degree $2(n - 1)$ such that
$$
  P(T_1, T_2, T_3, T_4) \, = \,  P(T_3, T_4, T_1, T_2) \, .
$$
$\GL^+ (2, \Rb)$ acts on $W_{n-1, n-1}$ by the tensor product of the
natural representations on the two factors $W_{n-1}$,
and the function $ F = F_{n-1, n-1} : \Hb \ra W_{n-1, n-1}$ is defined
by the formula
$$
 F_{n-1, n-1} (z) \, = \,  (zT_1 + T_2)^{n-1} \, (zT_3 + T_4)^{n-1} \, .
$$
For this specific choice, taking into account  \cite[Theorem 6.9]{Ste},
the statement of Theorem \ref{RadFi} becomes modified as follows.
\bigskip

\begin{theorem} \label{sv}
Let $m = 2n\geq 2$ and let ${\xb} \in \Qb^2/\Zb^2 $, with $\, {\xb} \neq 0$
if $\, m=2$.

$1^0$. For $\b =  \begin{pmatrix} a &b\\
 0 &d  \end{pmatrix} \in B (\Qb)$,
 \begin{equation} \nonumber
   \Phi^{(m)}_{\xb} (\b) \, = \, - \frac{ {\Bb}_m  (x_1)}{m} \,
   \int_{0}^{ \frac{b}{d} } (t T_1 + T_2)^{n-1} \, (tT_3 + T_4)^{n-1} \, dt  \, .
 \end{equation}

$2^0$. For $\s =  \begin{pmatrix} a &b\\
 c &d  \end{pmatrix} \in \SL(2, \Zb)$, with $c > 0$,
  \begin{equation} \nonumber
 \Phi^{(m)}_{\xb} (\s) = \left \{\begin{matrix} &-&
 \displaystyle{ \frac{ {\Bb}_m  (x_1)}{m} \,
   \int_{0}^{ \frac{a}{c} }} (t T_1 + T_2)^{n-1} \, (tT_3 + T_4)^{n-1} \,dt
   \qquad \qquad \qquad \qquad \quad  \cr \cr
&-&\displaystyle{\frac{ {\Bb}_m  (a x_1 + cx_2)}{m} \,
   \int_{-\frac{d}{c}}^{0 }} \left(t (aT_1 + cT_2) + bT_1 + d T_2 \right)^{n-1}
   \qquad \qquad \qquad \cr
  &\,& \qquad \qquad \qquad \qquad \cdot
   \left(t (aT_3 + cT_4) + bT_3 + d T_4 \right)^{n-1}  \,dt
   \qquad \qquad \cr \cr
&+& \displaystyle{\sum_{k=0}^{n-1} \sum_{\ell=0}^{n-1}} (-1)^{k + \ell}
 \begin{pmatrix} n-1\\
 k  \end{pmatrix}
  \begin{pmatrix} n-1\\
 \ell  \end{pmatrix}
  S^{(m-k- \ell -1, k+\ell +1)}_{\xb} (\frac{a}{c}) \qquad \qquad \cr
&\,& \qquad \cdot
 T_1^k \, (a T_1 + cT_2)^{n-k-1} \, T_3^{\ell} \, (a T_3 + cT_4)^{n-\ell-1} . \quad
\end{matrix} \right.
\end{equation}
\end{theorem}
\bigskip

\begin{remark}  {\rm Reformulating a result of Siegel \cite{Si1, Si2},
Stevens has shown \cite[\S 7]{Ste} that the Eisenstein cocycle
$\Phi$ of Theorem \ref{sv}
can be used to calculate the values at nonpositive integers of partial zeta
functions over a real quadratic field. This is achieved by specializing $\Phi$
to a certain Eisenstein series $E$, then
evaluating it at the element $\s \in \SL(2, \Zb)$
that represents the
action of a certain unit, and finally computing the polynomial  $\Phi_E (\s)$
on the basis elements and their conjugates.
It is intriguing to observe that the above transgressive construction
confers $\Phi$ the secondary status reminiscent of the Borel regulator invariants
that enter in the expression of the special values at non-critical points of
$L$-functions associated to number fields (\cite{KZ, Za1, Za2}).}
\end{remark}


\vspace{1cm}

\begin{thebibliography}{04}


\bibitem{As} Asai, T., The reciprocity of Dedekind sums and the factor set
for the universal covering group of $\SL(2, \Rb)$,
\textit{Nagoya Math. J.} \textbf{37} (1970), 67-80.

\bibitem{BG} Barge, J. and Ghys, E., Cocycles d'Euler et de Maslov,
\textit{Math. Ann.} \textbf{294} (1992), 235--265.

\bibitem{C1} Connes, A., Noncommutative differential geometry,
 \textit{Inst. Hautes Etudes Sci. Publ. Math.} \textbf{62} (1985), 257-360.

\bibitem{Cncg} Connes, A., \textbf{Noncommutative Geometry}, Academic
Press, 1994.

\bibitem{C2} Connes, A., Cyclic cohomology, quantum group symmetries
and the local index formula for $SU_q(2)$,
\textit{J. Inst. Math. Jussieu} \textbf{3} (2004), 17-68.


\bibitem{CMhti} Connes, A. and Moscovici, H., Hopf algebras, cyclic
cohomology and the transverse index theorem, \textit{Commun. Math.
Phys.} \textbf{198} (1998), 199-246.

\bibitem{CMcchs} Connes, A. and Moscovici, H.,  Cyclic
cohomology and Hopf algebra symmetry, \textit{Letters Math. Phys.}
\textbf{52}  (2000), 1--28.

\bibitem{CMmha} Connes, A. and Moscovici, H.,  Modular Hecke algebras
and their Hopf symmetry, \textit{Moscow Math. J.}
\textbf{4}  (2004).

\bibitem{CMrc} Connes, A. and Moscovici, H.,  Rankin-Cohen brackets
and the Hopf algebra of transverse geometry, \textit{Moscow Math. J.}
\textbf{4}  (2004), 111--130.

\bibitem{DHZ} Dupont, J., Hain, R., Zucker, S.,
Regulators and characteristic classes of flat bundles.
\textit{The arithmetic and geometry of algebraic cycles}
(Banff, AB, 1998),
47--92, CRM Proc. Lecture Notes, 24,
Amer. Math. Soc., Providence, RI, 2000.

\bibitem{Go} Gorokhovsky, A., Secondary classes and cyclic cohomology
of Hopf algebras,
\textit{Topology}, \textbf{41} (2002), 993 -- 1016.

\bibitem{HKRS} Hajac, P. M., Khalkhali, M., Rangipour, B., Sommerh\"auser, Y.,
Hopf-cyclic homology and cohomology with coefficients,
\textit{C. R. Math. Acad. Sci. Paris}, \textbf{338} (2004), 667--672.

\bibitem{HWZ} Hall, R. R., Wilson, J. C., Zagier, D.,
Reciprocity formulae for general Dedekind-Rademacher sums,
\textit{Acta Arith.}, \textbf{73} (1995), 389--396.

\bibitem{He} Hecke, E., Theorie der Eisensteinschen Reihen h\"oherer
Stufe und ihre Anwendung auf Funktionentheorie und Arithmetik,
\textit{Abh. Math. Sem. Univ. Hamburg} \textbf{5} (1927),
199-224\, ; also in \textbf{Mathematische Werke}, Third edition,
Vandenhoeck \& Ruprecht, G\"ottingen, 1983.

\bibitem{KR} Khalkhali, M., Rangipour, B., Cup products in
Hopf-cyclic cohomology,
\textit{C. R. Math. Acad. Sci. Paris}, \textbf{xxx} (200x), xxx--xxx.

\bibitem{KZ} Kontsevich M. and Zagier, D., Periods,
\textbf{Mathematics unlimited---2001 and beyond}, 771--808,
Springer, Berlin, 2001.


\bibitem{Ku} Kubota, T., Topological covering of $\SL(2)$ over a local
field, \textit{J. Math. Soc. Japan}, \textbf{19} (1967),
231-267.

\bibitem{Me1} Meyer, C., {\bf Die Berechnung der Klassenzahl Abelscher
Korper \"uber quadratischen Zahlkorpern}, Akademie-Verlag, Berlin,
1957.

\bibitem{Me2} Meyer, C., \"Uber einige Anwendungen
Dedekindsche Summen,
\textit{J. Reine Angew. Math.}, \textbf{198} (1957), 143-203.


\bibitem{Nak} Nakamura, H., Generalized Rademacher functions
and some congruence properties, \textbf{Galois Theory and Modular Forms},
Kluwer Academic Publishers, 2003, 375-394.

\bibitem{Pet} Petersson, H., Zur analytischen Theorie der Grenzkreisgruppen I,
 \textit{Math. Ann.}, \textbf{115} (1938), 23-67.

\bibitem{Ra} Rademacher, H., Zur Theorie der Modulfunktionen,
 \textit{J. Reine Angew. Math.}, \textbf{167} (1931), 312-366.

\bibitem{RG} Rademacher, H. and Grosswald, E., \textbf{Dedekind Sums},
The Carus Mathematical Monographs, No. \textbf{16}, Mathematical
Association of America, 1972.

\bibitem{Sc} Sczech, R., Eisenstein cocycles for $\GL_2(\Qb)$ and values of
$L$-functions in real quadratic fields,
\textit{Comment. Math. Helvetici}, \textbf{67} (1992), 363-382.

\bibitem{Si1} Siegel, C. L., Bernoullische Polynome und quadratische
Zahlk\"orper,
\textit{Nachr. Akad. Wiss. G\"ottingen Math.-physik}, \textbf{2} (1968), 7-38.

\bibitem{Si2} Siegel, C. L., \"Uber die Fourierschen Koeffizienten von
Modulformen,
\textit{Nachr. Akad. Wiss. G\"ottingen Math.-physik}, \textbf{3} (1970), 15-56.

\bibitem{St} Stevens, G., \textbf{Arithmetic on modular curves} ,
Progress in Mathematics, \textbf{20}, Birkh{\"a}user, Boston, MA,
1982.

\bibitem{Ste} Stevens, G., \textbf{The Eisenstein measure and real quadratic
fields}, in \textbf{Proceedings of the International Conference on Number Theory}
(Quebec, 1987), 887-927, de Gruyter, Berlin, 1987.

 \bibitem{Za1} Zagier, D., Valeurs des fonctions zeta des corps quadratiques
 r\'eels aux entiers n\'egatifs,
\textit{Ast\'erisque}, \textbf{41-42} (1977), 135-151.

\bibitem{Za2} Zagier, D., Polylogarithms, Dedekind zeta functions and the
algebraic $K$-theory of fields, In:
\textit{Arithmetic Algebraic Geometry}, Progr. Math.
\textbf{89}, Birkh\"auser (1991), 391-430.


\end{thebibliography}
\end{document}